\documentclass[12pt]{article}

\usepackage{amsmath}
\usepackage{amssymb}
\usepackage{amsthm}
\usepackage{a4wide}
\usepackage[english]{babel}

\numberwithin{equation}{section}

\newtheorem{Th}{Theorem}
\newtheorem{Lem}{Lemma}
\newtheorem{Prop}{Proposition}

\newcommand{\Dem}{\noindent{\bf Proof }}

\renewcommand{\P}{\mathbb{P}}
\newcommand{\N}{\mathbb{N}}
\newcommand{\Z}{\mathbb{Z}}
\newcommand{\Q}{\mathbb{Q}}
\newcommand{\R}{\mathbb{R}}
\newcommand{\C}{\mathbb{C}}

\newcommand{\etoile}{^\ast}

\newcommand{\Span}{{\rm Span}}

\newcommand{\Pti}{\widetilde P}
\newcommand{\Uti}{\widetilde U}

\newcommand{\eps}{\varepsilon}

\newcommand{\om}{\omega}
\newcommand{\eneq}{\end{equation}}
\newcommand{\lam}{\lambda}

\newcommand{\una}{\{1,\ldots,a\}}
\newcommand{\unN}{\{1,\ldots,N\}}
\newcommand{\zeronsN}{\{0,\ldots,n/N\}}
\newcommand{\zeroNmu}{\{0,\ldots,N-1\}}
\newcommand{\unq}{\{1,\ldots,q\}}

\newcommand{\calN}{\mathcal{N}}
\newcommand{\cut}{\Omega}
\newcommand{\hol}{{\mathcal H}}
\newcommand{\ord}{{\rm ord}}
\newcommand{\dd}{{\rm d}}
\newcommand{\rk}{{\rm rk}}

\newcommand{\Card}{\#}
\newcommand\tra{ \ ^t  }

\newcommand{\cstun}{c_1}
\newcommand{\cstnum}{c_2}
\newcommand{\csttr}{c_3}
\newcommand{\cstdesig}{c(\sigma)}

\newcommand{\dz}{d_0}

\newcommand{\combitiny}[2]{{\tiny \left( \! \!  \begin{array}{c} #1 \\ #2 \end{array} \! \! \right)}}
\newcommand{\combi}[2]{{  \left( \begin{array}{c} #1 \\ #2 \end{array} \right)}}

\newcommand{\calP}{{\mathcal P}}
\newcommand{\Asig}{{\mathcal A}_\sigma}
\newcommand{\Ksig}{{\mathcal K}_\sigma}
\newcommand{\Fsig}{F_\sigma}
\newcommand{\Gsig}{G_\sigma}

\newcommand{\kappasig}{\kappa_\sigma}
\newcommand{\rsig}{r_\sigma}

\newcommand{\Ainf}{{\mathcal A}_\infty}
\newcommand{\kappainf}{\kappa_\infty}
\newcommand{\rinf}{r_\infty}

\newcommand{\Azero}{{\mathcal A}_0}

\newcommand{\Fzero}{F_0}
\newcommand{\Gzero}{G_0}
\newcommand{\pizero}{\pi_0}
\newcommand{\unm}{\{1,\ldots,m\}}

\newcommand{\calY}{{\mathcal Y}}
\newcommand{\calYti}{\widetilde{\mathcal Y}}
\newcommand{\Mti}{\widetilde{M}}
\newcommand{\loga}{{\mathcal L}_\alpha}
\newcommand{\Deltati}{\widetilde \Delta}
\newcommand{\gammaun}{\gamma_1}
\newcommand{\gammade}{\gamma_2}
\newcommand{\gammatr}{\gamma_3}
\newcommand{\oun}{{\varpi_1}}
\newcommand{\ode}{{\varpi_2}}
\newcommand{\otr}{{\varpi_3}}
\newcommand{\evalpha}{{\rm ev}_\alpha}
\newcommand{\AAA}{E}

\newcommand{\calB}{{\mathcal{B}}}

\renewcommand{\Re}{{\rm Re \,}}

\title{Shidlovsky's multiplicity estimate and Irrationality of zeta values}
\author{St\'ephane Fischler}
\date{\today}

\begin{document}

\newcommand{\pb}[1]{{\bf #1}}

\newcommand{\Li}{{\rm Li}}
\newcommand{\norm}[1]{\lVert #1\rVert}

\newcommand{\ttt}{i_0}
\maketitle

\begin{abstract}
In this paper we follow the approach of Bertrand-Beukers (and of later work of Bertrand), based on differential Galois theory, to prove a very general version of Shidlovsky's lemma that applies to Pad\'e approximation problems at several points, both at functional and numerical levels (i.e., before and after evaluating at a specific point).  This allows us to obtain a new proof of the Ball-Rivoal theorem on irrationality of infinitely many values of Riemann zeta function at odd integers, inspired by the proof of the Siegel-Shidlovsky theorem on values of $E$-functions: Shidlovsky's lemma is used to replace Nesterenko's linear independence criterion   with Siegel's, so that no  lower bound is needed  on the linear forms in zeta values. The same strategy provides a new proof, and a refinement, of Nishimoto's theorem on values of $L$-functions of Dirichlet characters.
\end{abstract}

 \bigskip
 
\noindent MSC 2010 : 11J72 (Primary);  11M06, 34M03 (Secondary).

\section{Introduction}

The motivation of this paper comes from the Ball-Rivoal theorem: $\zeta(s)$ is irrational for infinitely many odd integers $s\geq 3$. Its proof is based on {\em explicit} Pad\'e approximation to polylogarithms. In order to try to generalize this result to  other  functions, it would be natural to use {\em non-explicit} Pad\'e approximation instead, for instance through Siegel's lemma. Several difficulties arise; the first one is the need for a lower bound on the linear forms in zeta values, in order to apply Nesterenko's linear independence criterion: such a lower bound cannot be obtained from a
 non-explicit construction. 

In this paper we overcome this difficulty by giving a new proof of the Ball-Rivoal theorem in which no lower bound on the linear forms is used. Indeed Nesterenko's linear independence criterion is replaced with Siegel's combined with a multiplicity estimate, namely a new generalization of Shidlovsky's lemma. We combine  an explicit construction of the linear forms with the strategy used to prove the Siegel-Shidlovsky theorem on values of $E$-functions (see for instance \cite[Chapter 3]{Shidlovsky}). 

\bigskip

Let $q$ be a positive integer, and $A \in M_q(\C(z))$. We fix  $P_1,\ldots,P_q\in\C[z]$ and $n \in\N=\{0,1,2,\ldots\}$ such that 
  $\deg P_i \leq n$ for any $i$.  Then with any solution $Y = \tra (y_1,\ldots,y_q)$ of the differential system $Y'=AY$   is  associated a remainder $R(Y)$ defined   by 
$$R(Y)(z)  = \sum_{i=1}^q P_i(z) y_i(z).$$

\medskip

Let $\Sigma$ be a finite subset of $  \C\cup\{\infty\}$. For each $\sigma\in\Sigma$, let 
$(Y_j)_{j\in J_\sigma}$ be a family of solutions of $Y'=AY$ such that the functions $R(Y_j)$, $j\in J_\sigma$, are $\C$-linearly independent and   holomorphic at $\sigma$; here $\sigma\in  \C\cup\{\infty\}$ might be a singularity of  the differential system $Y'=AY$. We agree that $J_\sigma=\emptyset$ if  $\sigma\not\in    \Sigma$, and   let $M(z) = [P_{k,i}(z)]_{1\leq i,k \leq q} \in M_q(\C(z))$ where the rational functions    $P_{k,i}\in\C(z)$ are defined for $k \geq 1$ and $1 \leq i \leq q$  by
\begin{equation} \label{eqdefpki}
\left( \begin{array}{c} P_{k,1} \\ \vdots \\ P_{k,q} \end{array}\right) = \left(\frac{\dd}{\dd z} + \tra A\right)^{k-1} 
 \left( \begin{array}{c} P_{1} \\ \vdots \\ P_{q} \end{array}\right).
 \eneq
 Obviously the poles of the coefficients $P_{k,i}$ of $M$ are among those of $A$.

\bigskip

The following multiplicity estimate appears essentially (see below) in \cite[Th\'eor\`eme 2]{DBShid}.

\begin{Th} \label{thz} There exists  a  positive constant  $\cstun$, which depends only on    $A$ and $\Sigma$,
 such that if 
\begin{equation} \label{eqhypdetnn}
\sum_{\sigma\in\Sigma}\sum_{j\in J_{\sigma}} \ord_\sigma(R(Y_j)) \geq (n+1)  q - n  \Card J_\infty  -\tau 
\eneq
with $0 \leq \tau \leq n - \cstun$, then  $\det M(z)$ is not identically zero.
\end{Th} 

The special case where $\Sigma = \{0\}$, $\Card J_0=1$, and $Y_j$ is analytic at 0 is essentially Shidlovsky's lemma (see \cite[Chapter 3, Lemma 8]{Shidlovski}). When $\Sigma\subset\C$, $\Card J_\sigma=1$ for any $\sigma$, and all functions $Y_j$ are obtained by analytic continuation from a single one, analytic at all $\sigma\in\Sigma$, this result was proved by Bertand-Beukers \cite{BB} with more details on the constant $\cstun$. 
Then Bertrand has allowed  \cite[Th\'eor\`eme 2]{DBShid} an arbitrary number of solutions at each $\sigma$, proving Theorem \ref{thz} under the additional assumptions that $\infty\not\in\Sigma$ and  the functions $Y_j$, $j \in J_\sigma$, are analytic at  $\sigma$.

Our proof of Theorem \ref{thz} (like that of  \cite[Th\'eor\`eme 2]{DBShid})  follows the strategy of \cite{BB}, based on differential Galois theory. 
The point is that we allow $\Sigma$ to contain $\infty$, and/or singularities of the  differential system $Y'=AY$: only the  remainders $R(Y_j)$ are assumed to be holomorphic at $\sigma$ (but not the functions $Y_j$, and not at points $\sigma'\in\Sigma$ distinct from $\sigma$). These features make  Theorem \ref{thz} general enough to cover essentially all Pad\'e approximation problems related to polylogarithms we have found in the literature, for instance the ones of Beukers \cite{BeukersLNM, BeukersBolyai}, Sorokin \cite{Sorokin94, Sorokinpi, SorokinApery}, and those of \cite{FR}.   In such a setting,   $\tau$  in Eq. \eqref{eqhypdetnn} appears as  the difference between the number of unknowns and the number of equations.

\bigskip

Then we evaluate at a point $\alpha$, going from functional to numerical linear forms (see \cite[Chapter 3, Lemma 10]{Shidlovski} for the classical setting). The point here is that we allow $\alpha$ to be   a singularity of the differential system $Y'=AY$, and/or an element of   $\Sigma$ (in our proof of the Ball-Rivoal theorem, $\alpha$ is both).

\begin{Th} \label{thzeronum}
There exists  a  positive constant   $\cstnum$, which depends only on    $A$ and $\Sigma$, with the following property. Assume that, for some $\alpha\in\C$:
  \begin{itemize}
\item[$(i)$] The differential system $Y'=AY$ has a basis of local solutions at $\alpha$ in \mbox{$\C[\log(z-\alpha)][[z-\alpha]]$.}
\item[$(ii)$] All rational functions  $P_{k,i}$, with $1\leq i \leq q$ and $1 \leq k \leq \tau + \cstnum$, are holomorphic at $z=\alpha$.
 \item[$(iii)$]    Eq. \eqref{eqhypdetnn} holds for some $\tau $ with $0 \leq \tau \leq n - \cstun$.
\end{itemize}
Then the matrix  $[P_{k,i}(\alpha)]_{1\leq i  \leq q, 1\leq k \leq \cstnum} \in M_{q,\cstnum}(\C)$ has rank at least $q-\Card J_\alpha$.
\end{Th}

If $\alpha$ is a singularity, assertion $(i)$ means it is regular and all exponents at $\alpha$ are integers. As far as we know, this result is the first general one in which $\alpha$ is allowed to be a singularity. The case where $\alpha$ is not a singularity is much easier, and assumptions $(i)$ and $(ii)$ are then  trivially satisfied.

 If $\alpha\not\in\Sigma$ then $J_\alpha=\emptyset$ so that we obtain a matrix of maximal rank $q$. On the opposite, if $\alpha\in\Sigma$  then $\Card J_\alpha$ linearly independent linear combinations of the rows of the matrix $[P_{k,i}(z)]_{i,k}$ are holomorphic at $\alpha$ and (probably) vanish at $\alpha$: the lower bound $q-\Card J_\alpha$ is best possible.

\bigskip

Using a zero estimate such as Theorem \ref{thzeronum} is the key point in the classical proof of the Siegel-Shidlovsky theorem on values of $E$-functions.  Following a different but similar strategy,  
Nikishin constructed explicitly \cite{Nikishin} linearly independent linear forms in 1, $\Li_1(\alpha)$, \ldots, $\Li_a(\alpha)$ to prove that these numbers are linearly independent over $\Q$ when $\alpha = u/v$ is a rational  number with $v$ sufficiently large in terms of $|u|$. His approach  was used by several authors, including Marcovecchio \cite{Marcovecchio} to bound from below the dimension of the $\Q$-vector space spanned by these numbers, for any fixed algebraic number $\alpha$ with $|\alpha| < 1$ (thereby generalizing to non-real numbers $\alpha$ Rivoal's result \cite{Rivoalpolylogs} based on Nesterenko's linear independence  criterion). The zero estimate used by Marcovecchio is similar to   Theorem \ref{thzeronum} but deals only with a specific situation in which (essentially) $\tau=1$ in Eq. \eqref{eqhypdetnn}, $\alpha\not\in\Sigma$, and $\alpha$ is not a singularity. Moreover he does not define $P_{k,i}$ for $k\geq 2$ using  Eq. \eqref{eqdefpki} (i.e., differentiating the linear forms as in  the  proof of the Siegel-Shidlovsky theorem): following Nikishin he uses an additional parameter instead.

\bigskip

In this paper we use Theorem \ref{thzeronum} to obtain a new proof, and a refinement,   of the following result of Nishimoto   \cite{Nishimoto} on 
 $L$-functions $L(\chi,s) = \sum_{n=1}^\infty \frac{\chi(n)}{n^s}$ associated with Dirichlet characters $\chi$. He proved it with $d$ instead of $N$ in the lower bound \eqref{eqminocaract}; see \S \ref{subsec31nv} for this easy improvement.

\begin{Th} \label{thminocaract}
Let  $\chi$ be  a Dirichlet character modulo $d$, of conductor $N$. Let $p\in\{0,1\}$ and $a\geq 2$. Denote by $\delta_{\chi,p,a}$ the dimension of the $\Q$-vector space spanned by 1 and the numbers $L(\chi,s) $ with $2\leq s \leq a$ and $s\equiv p \bmod 2$. Then
\begin{equation} \label{eqminocaract}
 \delta_{\chi,p,a} \geq \frac{1+o(1)}{N+\log 2} \log a
 \eneq
 where $o(1)$ is a sequence that depends on $N$ and $a$, and tends to 0 as $a\to\infty$ (for any $N$).
 \end{Th}

If $p$ and $\chi$ have the same parity then $L(\chi,s) \pi^{-s}$ is a  non-zero algebraic number for any $s\geq 2$ such that $s\equiv p \bmod 2$ (see for instance \cite[Chapter VII, \S 2]{Neukirch}):  this result is interesting when $p$ and $\chi$ have opposite parities.

Nishimoto's proof is similar to Ball-Rivoal's, except that obtaining the lower bound necessary to apply Nesterenko's  criterion is very technical: the saddle point method has to be used because cancellations take place (see \cite{Nash}).  In this paper we present an alternative proof of Theorem \ref{thminocaract}, based on the zero estimate stated above. It makes it unnecessary to use the saddle point method, since Siegel's criterion is applied instead of Nesterenko's. In the special case $d=N=1$ (so that $\chi(n) =1$ for any $n$, and $L(\chi,s) = \zeta(s)$) this is exactly the proof of the Ball-Rivoal theorem mentioned above.

\bigskip

We also obtain the following refinement of Theorem \ref{thminocaract}, by improving the arithmetic estimates.

\begin{Th} \label{thminocaractnv}
In the setting of Theorem \ref{thminocaract}, if $N$ is a multiple of 4 then Eq. \eqref{eqminocaract} can be replaced with
$$
 \delta_{\chi,p,a} \geq \frac{1+o(1)}{(N/2)+\log 2} \log a.
 $$
 \end{Th}
When $\chi$ is the non-principal character mod $d=4$, so that $N=4$, this result was proved by Rivoal-Zudilin \cite{Catalan} as a first step towards the (conjectural) irrationality of Catalan's constant $L(\chi,2) = \sum_{k=0}^\infty \frac{(-1)^k}{(2k+1)^2}$.  

\bigskip

The structure of this paper is as follows. We first sketch in \S \ref{secBR} our  proof of the Ball-Rivoal theorem. Then \S \ref{sec2} is devoted to Shidlovsky's lemma: we prove Theorems \ref{thz} and \ref{thzeronum}. At last, in \S \ref{secdio} we prove in details a general result which contains Theorem \ref{thminocaract}, Theorem \ref{thminocaractnv}, and the Ball-Rivoal theorem.

\section{A new proof of the Ball-Rivoal theorem} \label{secBR}

We sketch in this section the new proof of the Ball-Rivoal theorem obtained as a special case of the proof of Theorem \ref{thprinc} in \S \ref{secdio} below (namely $N=1$, $f(n)=1$ for any $n$, $p=1$, $z_0=1$, $\ttt=2$, $\xi_1=0$, and $\xi_j = \zeta(j)$ for any $j\geq 2$). Of course we refer to \S \ref{secdio} for more details.

\medskip

Let $a$, $r$, $r'$, $n$ be such that $a$ is odd and $r,r' < a/2$. It turns out that the best estimates come from the case where $r$ and $r'$ have essentially the same size, so we shall restrict in \S \ref{secdio} to the case $r'=r$; however the proof works in the same way if $r'\neq r$. Consider the rational function
$$F(t) = n!^{a-r-r'} \frac{(t-rn)_{rn} (t+n+1)_{r'n}}{(t)_{n+1}^a}$$
where $(\alpha)_k = \alpha ( \alpha+1)\ldots (\alpha+k-1)$ is Pochhammer's symbol, and let 
$$S_0(z) = \sum_{t=n+1}^\infty F(-t) z^t, \hspace{2cm} 
S_\infty(z) = \sum_{t=1}^\infty F( t) z^{-t} .
$$
For any $k\geq 1$ we let 
\begin{equation} \label{eqlamBR}
\Lambda_k = S_0^{(k-1)}(1) - S_\infty^{(k-1)}(1) ,
\eneq
where $S^{(k-1)}$ is the $(k-1)$-th derivative of $S$.
We shall use a symmetry phenomenon to get rid of even zeta values, but it does not appear exactly as in the original proof of Ball-Rivoal.  Indeed, even if $r'=r$, $S_0^{(k-1)}(1) $ and $ S_\infty^{(k-1)}(1)$ involve both   odd and even zeta values when $k\geq 2$: they are values at $z=1$ of hypergeometric series which are no more well-poised. The cancellation of even zeta values comes at a different stage, by considering $\Lambda_k$ in Eq. \eqref{eqlamBR}. Indeed there exist integers $s_{k,i}$, $2\leq i \leq a$, and $u_k$, $v_k$ such that for any $k\leq (a-r-r')n+a-1$, we have both
$$d_n^a S_0^{(k-1)}(1)  = u_k + \sum_{i=2}^a (-1)^i s_{k,i} \zeta(i)$$
and
$$d_n^a S_\infty^{(k-1)}(1)  = v_k + \sum_{i=2}^a   s_{k,i} \zeta(i)$$
where $d_n = {\rm lcm}(1,2,\ldots,n)$, so that 
 $d_n^a\Lambda_k  =d_n^a S_0^{(k-1)}(1) - d_n^a S_\infty^{(k-1)}(1) $ is a $\Z$-linear combination of 1 and odd zeta values:
 $$ d_n^a\Lambda_k  = s_{k,a+1} - 2  \sum_{2 \leq i \leq a \atop i \mbox{ {\tiny odd}} }  s_{k,i} \zeta(i),$$
  with $s_{k,a+1} = u_k-v_k$. Using Theorem \ref{thzeronum} we prove that the matrix $[s_{k,i}]_{2\leq i \leq a+1, 1\leq k \leq \cstnum}$ has maximal rank, equal to $a$ (see below). This enables one to apply Siegel's linear independence criterion (see \S \ref{subsec36}) instead of Nesterenko's: no lower bound on $|\Lambda_k|$ is needed. The upper bounds  on $|s_{k,j}|$ and  $|\Lambda_k|$ are essentially the same as in the proof of Ball-Rivoal, so that we obtain the same lower bound:
  $$\dim_\Q \Span_\Q (1, \zeta(3), \zeta(5), \ldots, \zeta(a)) \geq \frac{1+o(1)}{1+\log 2} \log a.$$
  
\bigskip

Let us focus now on the functional aspects of this proof, which play an important role (whereas the proof of Ball-Rivoal  can be written with $z=1$ throughout). For simplicity we restrict ourselves to the case $r'=r$.  The functions $S_0(z)$ and $S_{\infty}(z)$ are solutions of the following Pad\'e approximation problem: find polynomials $P_1$, \ldots, $P_{a+2}$ of degree at most $n$ such that:
\begin{equation} \label{eqpadeFR}
\left\{
\begin{array}{rcl}
S_0(z) := P_{a+1} (z) +   \sum_{i=1}^a P_i(z) (-1)^i \Li_i(  z) = O(z^{(r+1)n+1}), \hspace{0.5cm} z\to 0,\\
S_{\infty}(z) := P_{a+2} (z) +   \sum_{i=1}^a P_i(z)   \Li_i(1/z) = O(z^{ -rn-1}), \hspace{0.5cm} z\to \infty, \\
  \sum_{i=1}^a P_i(z)  (-1)^{i-1}\frac{(\log z)^{i-1}}{(i-1)!} = O((z-1)^{(a-2r)n +a-1}),  \hspace{0.5cm} z\to 1.
  \end{array}
\right.
\end{equation}
This is exactly the  Pad\'e approximation problem of \cite[Th\'eor\`eme 1]{FR}: it has a unique solution up to proportionality, $(n+1)(a+2)$ unknowns and $(n+1)(a+2)-1$ equations. Let $A \in M_{a+2}(\C(z))$ denote the following matrix:
$$A = 
\left[
\begin{matrix}
0 & 0 & 0 & \ldots & 0 & 0 & \frac1{z-1} & \frac{1}{z(1-z) }\\
\frac{-1}z  & 0  & 0 & \ldots  & 0  & 0  & 0  & 0\\
0  &   \frac{-1}z  & 0 & \ldots  & 0  & 0  & 0 & 0 \\
 0  & 0  &  \frac{-1}z  &  \ldots& 0  & 0  & 0  & 0\\
\vdots & \vdots & \vdots & \ddots & \vdots & \vdots & \vdots & \vdots \\
0 & 0 & 0 & \ldots &  \frac{-1}z  &  0  & 0  & 0\\
0 & 0 & 0 & \ldots &  0  &  0  & 0  & 0\\
0 & 0 & 0 & \ldots &  0  &  0  & 0  & 0
\end{matrix}
\right]
$$
and consider the following solutions of the differential system $Y'=AY$:
$$Y_0(z) = \tra ( -\Li_1(z), \, \Li_2(z), \, \ldots, (-1)^a\Li_a(z), 1, 0),$$
$$Y_\infty(z) = \tra (  \Li_1(1/z), \, \Li_2(1/z), \, \ldots,  \Li_a(1/z), 0, 1),$$
$$Y_1(z) = \tra (1,  -\log z , \, \frac{(\log z)^2}2 , \, \ldots, (-1)^{a-1}  \frac{(\log z)^{a-1}}{(a-1)!},  0,  0).$$
Let $\Sigma = \{0,1,\infty\}$ and $J_0=\{0\}$, $J_1=\{1\}$, 
$J_\infty=\{\infty\}$. Then with the notation of the introduction, we have  $R(Y_0) = S_0(z)$,  $R(Y_\infty) = S_\infty(z)$,  and $R(Y_1)$ is the left hand side of the third equation of \eqref{eqpadeFR}; Eq. \eqref{eqhypdetnn} stated in the introduction holds with $\tau=1$ as a consequence of the Pad\'e approximation problem \eqref{eqpadeFR}. In general, $\tau$ corresponds in Eq.~\eqref{eqhypdetnn} to the difference between the number of unknowns and the number of equations. To apply Theorem \ref{thzeronum} it is not useful to prove that the problem has a unique solution up to proportionality: the upper bound $\tau \leq n/2$, for instance, would be  sufficient since $n$ is taken arbitrarily large. 

Defining $P_{k,i}$ as in the introduction by Eq. \eqref{eqdefpki}, it is well-known (see \cite[Chapter 3, \S 4]{Shidlovski})¤ that for any $k\geq 1$, 
\begin{equation} \label{eqokenunBR} 
S_{0 }^{(k-1)}(z) = P_{k,a+1}(z) +\sum_{i=1}^a P_{k,i} (z) (-1)^i \Li_i(  z) \mbox{  and }
 S_{\infty }^{(k-1)} (z) = P_{k,a+2}(z) +\sum_{i=1}^a P_{k,i} (z)   \Li_i(1/z ). 
 \eneq
 Moreover $P_{k,i}$ is a rational function of which $0$ is the only possible pole if $i \leq a$. If $i=a+1$ or $i = a+2$, both 0 and 1 may be poles of $P_{k,i}$; but if $k \leq (a-2r )n+a-1$, the functions $S_{0 }^{(k-1)}(z) $ and $ S_{\infty }^{(k-1)} (z) $ have finite limits as $z\to1$ so that 1 is not a pole. 
 
 Finally Theorem \ref{thzeronum} applies at $\alpha=1$: the matrix $[P_{k,i}(1)]_{1\leq i \leq a+2, 1\leq k \leq \cstnum}$ has rank at least $a+1$. Actually $P_{k,1}(1)=0$ for any $k\leq (a-2r)n+a-1$ (which can be seen by letting $z$ tend to 1 in Eq. \eqref{eqokenunBR}) so that the first row of this matrix is zero (provided $n$ is large enough) and its rank is exactly $a+1$. Since the coefficients $s_{k,i}$ defined above are given by $s_{k,i} = d_n^a P_{k,i}(1)$ for $2\leq i \leq a$ and $s_{k,a+1} = d_n^a ( P_{k,a+1}(1) - P_{k,a+2}(1) )$, the matrix $[s_{k,i}]_{2\leq i\leq a+1, 1\leq k \leq \cstnum}$ has rank $a$: Siegel's criterion (stated and proved in \S \ref{subsec36}) applies.

\section{Zero estimates} \label{sec2}

In this section we prove Theorems \ref{thz} and \ref{thzeronum}. We start with the functional part of the proof (\S \ref{subsec21}), in which we follow the approach of Bertrand-Beukers \cite{BB} to generalize Shidlovsky's lemma (see Theorem \ref{thzerofct}). Then we deduce in \S \ref{subsec22} Theorems \ref{thz} and \ref{thzeronum} stated in the introduction: the important point is to evaluate at $\alpha$ which may be a singularity and/or an element of $\Sigma$.

\subsection{Functional zero estimate} \label{subsec21}

Throughout this section we consider a positive integer $q$ and a matrix $A \in M_q(\C(X))$. We let $P_1,\ldots,P_q\in\C[X]$ with $\deg P_i \leq n$ for any $i$. We also denote by $\cut$ a simply connected open subset of $\C$ in which $A$ has no pole. We assume that $\cut $ is obtained from $\C$ by removing finitely many half-lines, so that $\cut $ is dense in $\C$, and   denote  by $\hol$ the space of functions holomorphic on $\cut$. A solution $Y$ of the differential system $Y'=AY$ will always be a column matrix in $M_{q,1}(\hol)$, identified with the corresponding element $(y_1,\ldots,y_q)$ of $\hol^q$. Since $P_1,\ldots,P_q$ are fixed, to such a solution is  associated a remainder $R(Y)$ defined on $\cut$ by 
$$R(Y)(z)  = \sum_{i=1}^q P_i(z) y_i(z).$$

Let $\Sigma$ be a finite subset of $\P^1(\C) = \C\cup\{\infty\}$. For each $\sigma\in\Sigma$, let 
$(Y_j)_{j\in J_\sigma}$ be a family of solutions of $Y'=AY$  such that:
\begin{itemize}
\item For any $j\in J_\sigma$, the function $R(Y_j)$ is holomorphic at $\sigma$.
\item The functions $R(Y_j)$, for $j\in J_\sigma$, are linearly independent over $\C$.
\end{itemize}
Here we do not assume that $\sigma \in \cut$: in the case  $\sigma\not\in\cut$ (for instance if $\sigma=\infty$), by {\em   $R(Y_j)$ is holomorphic at $\sigma$} we mean that $R(Y_j)$ can be continued analytically to a function  holomorphic at $\sigma$. Moreover, we denote by $\ord_\sigma(R(Y_j))$ its order of vanishing at $z=\sigma$.

The point is that we do not assume any relation (or lack of relation) between the families $(Y_j)_{j\in J_\sigma}$ at distinct points $\sigma$, except of course that all are solutions of the same differential system.

At last, we let $J_\sigma=\emptyset$ when $\sigma\not\in\Sigma$.

\bigskip

Defining $M(z)$ and $P_{k,i}$ as   in the introduction, our functional multiplicity estimate is the following  generalization of Bertrand-Beukers' version of Shidlovsky's lemma; if $\infty\not\in\Sigma$  and  the functions $Y_j$, $j \in J_\sigma$, are analytic at  $\sigma$ it is   due to Bertrand  \cite[Th\'eor\`eme 2]{DBShid}. The constant $\cstun$ is the same as in Theorem~\ref{thz} (that we shall deduce from Theorem \ref{thzerofct} at the beginning of \S \ref{subsec22}).

\begin{Th} \label{thzerofct}
Let $\mu$ denote the order of  a non-zero differential operator  $L \in \C(z)[\frac{\dd}{\dd z}]$  such that $L(R(Y_j))= 0 $ for any $\sigma$ and any $j\in J_\sigma$. Then 
\begin{equation} \label{eq39}
\sum_{\sigma\in\Sigma}\sum_{j\in J_{\sigma}} \ord_\sigma(R(Y_j)) \leq (n+1) (\mu - \Card J_\infty) + \cstun
\eneq
where $\cstun$ is a constant that depends only on  $A$ and $\Sigma$.
\end{Th}

In the special case where $\Sigma\subset\C$, $J_\sigma$ consists of a single element $j_\sigma$, and the function $Y_{j_\sigma}$ is the same for all $\sigma$,   this is exactly \cite[Th\'eor\`eme 2]{BB} except that we did not try to make the constant $\cstun$ explicit
(we refer to \cite{BB}, and to \cite[Appendix of Chapter III]{Andre} in the Fuchsian case, for discussions on effectivity which are not relevant to our purposes). Indeed we have fixed a simply connected open subset $\Omega$ only for convenience: analytic continuation from a point of $\Sigma$ to another could be performed along any fixed path. 

\bigskip

Let us prove Theorem \ref{thzerofct} now, following the strategy of \cite{BB}.

Given $\sigma\in\P^1(\C)$, we let $\Asig$  denote the set of all finite sums
\begin{equation} \label{eqasig}
\sum_{\alpha\in E} \sum_{Q\in\calP} \sum_{j=0}^J u_{\alpha,Q,j} (z-\sigma)  (z-\sigma)^\alpha (\log (z-\sigma))^j \exp(Q( (z-\sigma)^{-1/q!}))
\eneq
where $E\subset \C$ and $\calP \subset\C[X]$ are finite subsets, $J\geq 0$, and $u_{\alpha,Q,j} (z-\sigma) \in \C [[  (z-\sigma)^{ 1/q!} ]]$ for any $\alpha$, $Q$, $j$. Here and below, we agree that $z-\sigma$ stands for $1/z$ si $\sigma=\infty$. Then the differential system $Y'=AY$ has a complete system of formal solutions in $\Asig^q$. Moreover  we let $\Ksig$  denote the fraction field of $\Asig$, and $\Fsig $ denote the differential subfield of $\Ksig$ generated over $\C(z)$ by all components of all solutions of $Y'=AY$  in $\Ksig^q$. Then the differential extension $\Fsig/ \C(z)$   is Picard-Vessiot, and we denote by $\Gsig$ its group of differential automorphisms.

\bigskip

To prove Theorem \ref{thzerofct} we may assume that $0\in\Sigma$, that $\mu$ is the minimal order of a non-zero  differential operator that annihilates  $R(Y_j)$  for any $j\in J_\sigma$ and any $\sigma\in\Sigma$, and that the coefficient of $(\frac{\dd}{\dd z})^\mu$ in $L$ is 1. 

Given $\sigma\in\Sigma$ and $j\in J_\sigma$, all components of $Y_j$ are holomorphic on the cut plane $\Omega$, and can be seen as elements of $\Azero$. Indeed, if 0 is a regular singularity (or an ordinary point) of the system $Y'=AY$ then all components of $Y_j$ have a generalized Taylor expansion at the origin in $\Azero$ (of the form \eqref{eqasig} with $\calP = \{0\}$). In the general case, we identify each   component  of $Y_j$ with its asymptotic expansion at 0 in a fixed large sector (see \cite{Ramis}). By definition of $\Fzero$, all   components of $Y_j$  (seen in  $\Azero$) belong to $\Fzero$ so that $R(Y_j) \in F_0$. We consider the $\C$-vector space $V \subset \Fzero$ spanned by the images $\gamma(R(Y_j))$ of all $R(Y_j)$, $j\in J_\sigma$, $\sigma\in\Sigma$,  under all $\gamma\in\Gzero$. Since the kernel of $L : \Fzero\to\Fzero $ is stable under $\Gzero$, we have $V\subset\ker L$ so that $m \leq \mu$, where $m  = \dim_\C V$. Let $(R_1,\ldots,R_m)$ be a basis of $V$, such that $R_i = \gamma_i ( R(Y_{j_i}))$ with $\gamma_i \in \Gzero$ and $j_i \in J_{\sigma_i}$ for any $i \in \unm$.

Arguing as in the proof of \cite[Proposition 3]{BB}, we have
$$Ly = \frac1{W(R_1,\ldots,R_m)} \, \det \left[ \begin{matrix} y & y' & \ldots & y^{(m)} \\
R_1 & R'_1 & \ldots & R^{(m)}_1 \\
\vdots & \vdots &  & \vdots &\\
R_m & R'_m & \ldots & R^{(m)}_m  
\end{matrix}\right]$$
where $W(R_1,\ldots,R_m)  = \det [R_i^{(j-1)}]_{1\leq i, j\leq m}$ is the wronskian determinant. In particular, we have $m = \mu  $ and $V = \ker L$. 

\bigskip

Now we claim that for any $\sigma\in \P^1(\C)$ there exist $\mu$ solutions $Y^{[\sigma,j]} = (y_1^{[\sigma,j]}, \ldots, y_q^{[\sigma,j]})$ of $Y'=AY$ in $\Fsig^q$, with $1\leq j \leq \mu$, such that $R(Y^{[\sigma,1]})$, \ldots, $R(Y^{[\sigma,\mu]})$ span the $\C$-vector space of solutions of $Ly=0$ in $\Fsig$. 
Indeed, as   in \cite[Corollaire]{BB}, using a differential isomorphism $\Fzero\to\Fsig$ we may assume $\sigma=0$. Then for any $i \in \unm$, $\gamma_i(Y_{j_i})\in\Fzero^q$ is a solution of $Y'=AY$ and $\pizero( \gamma_i(Y_{j_i})) = R( \gamma_i(Y_{j_i})) = \gamma_i(R( Y_{j_i})) = R_i$ so that  the claim is proved since $(R_1,\ldots,R_m)$ is a basis of $V = \ker L$.

\bigskip

Let us recall the following terminology from \cite{BB}: an element of $\Asig$ has rank $ \leq  \kappa \in \frac1{q!}\N$ and generalized order $\geq r$ if it is of the form \eqref{eqasig} with $\deg Q \leq q! \kappa$ for any $Q\in\calP$ and $\Re \alpha \geq r $ for any $\alpha\in E$. The differential operator $L$ has rank $ \leq  \kappa  $  at $\sigma$ and $(r_1,\ldots,r_\mu)\in\R^\mu$ is an admissible system of exponents of $L$ at $\sigma$ if the differential equation $Ly=0$ has a complete system of solutions $(y_1,\ldots,y_\mu)$ in $\Asig^\mu$ such that each $y_i$ has  rank $ \leq  \kappa $ and generalized order $\geq r_i$.

Given $\sigma\in\P^1(\C)$ all functions $y_i^{[\sigma,j]}$ with $1\leq i \leq q$ and $1 \leq j \leq \mu$   have rank $\leq \kappasig$ and generalized order $\geq \rsig$ for some $\kappasig \in \frac1{q!}\N$ and $\rsig\in\R$ which depend only on $A$ and $\sigma$ (see \cite[Proposition 1]{BB}). If $\sigma\neq\infty$, $R(Y^{[\sigma,j]}) = \sum_{i=1}^q P_i(z) y_i^{[\sigma,j]}(z) $ has rank $\leq \kappasig$ and generalized order $\geq \rsig$; these functions make up a complete system of solutions of $L$ (using the claim above). Moreover, if $\sigma\in\Sigma\setminus\{\infty\}$ then for any $j\in J_\sigma$ the function $R(Y_j)$ is holomorphic at $\sigma$,  so it can be seen as an element of $\Asig$ with rank $\leq 0$ and generalized order $\geq  \ord_\sigma(R(Y_j))$. Combining these $\C$-linearly independent solutions of $Ly=0$ with suitable functions $R(Y^{[\sigma,j]})$, 
  we obtain that $L$ has rank $\leq \kappasig $ at $\sigma$ and an admissible system of exponents of $L$ at $\sigma$ consists in $\rsig$ repeated $\mu-\Card J_\sigma$ times, and $ \ord_\sigma(R(Y_j))$ for each $j\in J_\sigma$. 
    In the same way, at infinity, for any $j\in J_\infty$ the function $R(Y_j)\in\Ainf$ has rank $\leq 0$ and generalized order $\geq  \ord_\infty(R(Y_j))$. To obtain a complete system of solutions of $Ly=0$ in $\Ainf$ we use also $\mu-\Card J_\infty$ functions $R(Y^{[\infty,j]})$, which have rank $\leq \kappasig$ and generalized order $\geq \rinf - n$ since $\deg P_i(z)\leq n $ for any $i \in \unq$. Therefore $L$ has  rank $\leq \kappainf $ at $\infty$ and an admissible system of exponents of $L$ at $\infty$ consists in $\rinf - n $ repeated $\mu-\Card J_\infty$ times, and $ \ord_\infty(R(Y_j))$ for each $j\in J_\infty$. 

So far we have found an upper bound on the rank of $L$, and an admissible system of exponents of $L$, at any $\sigma\in\Sigma$. Enlarging $\Sigma$ if necessary, we may assume that it contains $\infty$ and all poles of $A$. Then for any $\sigma\in \P^1(\C)\setminus\Sigma$ the differential system $Y'=AY$ has a complete system of solutions holomorphic at $\sigma$, and therefore the same property holds for the differential equation $Ly=0$ using the claim above. Accordingly $\Sigma$ contains $\infty$ and all non-apparent singularities of $L$, so that the Corollary of \cite[Th\'eor\`eme 3]{BB} provides an inequality involving upper bounds on the ranks of $L$  and   admissible systems of exponents of $L$  at  all points of $\Sigma$, namely:
\begin{eqnarray*}
&&(\mu-\Card J_\infty)(\rinf-n) + \Big( \sum_{j\in J_\infty} \ord_\infty (R(Y_j))\Big) - (\kappainf + 1)\mu(\mu-1)/2 \\
&&+ \sum_{\sigma\in\Sigma\setminus\{\infty\}} \Big[ (\mu-\Card J_\sigma) \rsig + \Big( \sum_{j\in J_\sigma} \ord_\sigma (R(Y_j))\Big) - (\kappasig + 1)\mu(\mu-1)/2 \Big] \leq -\mu(\mu-1)
\end{eqnarray*}
so that
$$ \Big( \sum_{\sigma\in\Sigma} \sum_{j\in J_\sigma} \ord_\sigma (R(Y_j)) \Big) - (n+1)( \mu-\Card J_\infty) \leq  \cstun$$
where $\cstun$ is a constant that can be written down explicitly in terms of $\Sigma$, $\mu$, $\kappasig$, $\rsig$ and $\Card J_\sigma $ for $\sigma\in\Sigma$. This concludes the proof of Theorem \ref{thzerofct}.

\subsection{Numerical zero estimate} \label{subsec22}

In this section we prove Theorems \ref{thz} and  \ref{thzeronum} stated in the introduction. The proof falls into 3 steps; the first one is Theorem \ref{thz}.

\bigskip

\noindent {\bf Step 1: } $M(z) \in M_q(\C(z))$ is an invertible matrix.

\medskip

As in \cite{Shidlovski},  if $M$ is singular in $M_q(\C(z))$ then there is a non-trivial linear relation with coefficients in $\C(z)$ between the $\rk(M)+1$ first columns of $M$; this provides a differential operator $L$ of order $\mu = \rk(M)$ to which Theorem \ref{thzerofct} applies, in contradiction
 with Eq.~\eqref{eqhypdetnn} since $\tau \leq n -\cstun$.
   Indeed, for any solution $Y$ of the differential system $Y'=AY$ we have
 $$\tra Y M = \left[ \begin{array}{cccc} R(Y) & R(Y)' & \ldots & R(Y)^{(q-1)}\end{array}\right].$$

\bigskip

\noindent {\bf Step 2: } Determination of $\det M(z)$ up to factors of bounded degree.

\medskip

\newcommand{\matsig}{T_\sigma}
\newcommand{\Nsig}{P_\sigma}
\newcommand{\Ninf}{P_\infty}

Let $S$ denote the set of finite singularities of the differential system $Y'=AY$, i.e. poles of coefficients of $A$. For any $s\in S$, let $N_s$ denote the maximal order of $s$ as a pole of a coefficient of $A$; let $N_s = 0$ for $s\in\C\setminus S$. Then Eq. \eqref{eqdefpki} shows that $(z-s)^{(k-1)N_s}P_{k,i}(z) $ is holomorphic at $z=s$ for any $k\geq 1$ and any $i\in\unq$. Therefore $\det M(z) \cdot \prod_{s\in S} (z-s)^{q(q-1)N_s}$ has no pole: is it a polynomial.

Now let $\sigma\in \Sigma$, and denote by $\matsig \in M_{\Card J_\sigma, q}(\hol)$ the matrix with rows $\tra Y_j$, $j \in J_\sigma$. 
The vector-valued functions $Y_j$, $j\in J_\sigma$, are linearly independent over $\C$ because the functions $R(Y_j)$ are; therefore 
they are the $ \Card J_\sigma$ first elements of a basis of solutions $\calB$ of the differential system $Y'=AY$. The wronskian determinant of $\calB$ may vanish at $\sigma$ if $\sigma$ is a singularity, but even in this case it cannot have generalized order $\geq c_0(\sigma)$ at $\sigma$ (with the terminology of \S \ref{subsec21}) where $c_0(\sigma)$ is a constant depending only on $A$ and $\sigma$ (not on $\calB$). On the other hand, all components of all elements of $\calB$ have  generalized order $\geq r_\sigma $ at $\sigma$ (as in  \S \ref{subsec21}). Therefore 
 there exists a subset $I_\sigma$ of $\unq$, with $\Card I_\sigma  = q- \Card J_\sigma$, such that the determinant of the submatrix of $\matsig$ corresponding to the columns indexed by $\unq \setminus I_\sigma$
 cannot have  generalized order $\geq c(\sigma)$ at $\sigma$,
 where $c(\sigma) = c_0(\sigma) - r_\sigma \Card I_\sigma$ depends only on $A$ and $\sigma$.

 Let $\Nsig \in M_q(\hol) $ denote the matrix of which the $\Card J_\sigma$ first rows are that of $\matsig$, and the other rows are the $\tra e_i$, $i\in I_\sigma$, where $(e_1,\ldots,e_q)$ is the canonical basis of $M_{q,1}(\C)$. Then   $\Nsig M$ has its first rows equal to $ \left[ \begin{array}{cccc} R(Y_j) & R(Y_j)' & \ldots & R(Y_j)^{(q-1)}\end{array}\right]$ with $j\in J_\sigma$, and its last rows equal to  $ \left[ \begin{array}{ccc} P_{1,i} \ldots P_{q,i}\end{array}\right]$ with $i \in I_\sigma$. Therefore all coefficients in the row corresponding to $j \in J_\sigma$ vanish at $\sigma$ with order at least $\ord_\sigma R(Y_j)-q+1$, and (if $\sigma\neq\infty$) all  coefficients in the row corresponding to $i \in I_\sigma$ are either holomorphic at $\sigma$, or have a pole of order at most $(q-1)N_\sigma$ is $\sigma\in S$. Since $N_\sigma = 0$ if $\sigma\not\in S$, we have for any $\sigma\in\Sigma\setminus\{\infty\}$:
$$\ord_\sigma \det(\Nsig M) \geq \Big( \sum_{j\in J_\sigma}\ord_\sigma R(Y_j)\Big) - (q-1)\Card J_\sigma -(q-1)N_\sigma (q-\Card J_\sigma).$$
Since $\det \Nsig$  cannot have  generalized order $\geq c(\sigma)$ at $\sigma$,  we obtain
$$\ord_\sigma \det( M) \geq \Big( \sum_{j\in J_\sigma}\ord_\sigma R(Y_j)\Big) - \cstdesig . $$

Now let
$$Q_2(z) = \Big( \prod_{s\in S} (z-s)^{q(q-1)N_s} \Big) \cdot \Big( \prod_{\sigma\in\Sigma\setminus\{\infty\}}(z-\sigma)^{\cstdesig}\Big)$$
so that $Q_2(z) \det M(z)$ is a polynomial and vanishes at any $\sigma \in\Sigma\setminus\{\infty\}$ with order at least $ \sum_{j\in J_\sigma}\ord_\sigma R(Y_j)$. To bound from above the degree of this polynomial, we define $\Ninf$ as above if $\infty\in\Sigma$, and let $\Ninf $ denote the identity matrix (and $J_\infty = \emptyset$) otherwise. Then we have $R(Y_j)^{(k-1)} = O(z^{-\ord_\infty R(Y_j)})$ as $|z|\to\infty$ for any $j\in J_\infty$ and any $k\geq 1$, and $P_{k,i}(z) = O(z^{n+(q-1)d})$ for any $i\in I_\infty$ and any $k \in \unq$ (where $d$ is greater than or equal to the degree of all coefficients of $A$). Therefore we have $\det M(z)  = O(z^u)$ as $|z| \to \infty$, with
$$u = (q-\Card J_\infty)(n+(q-1)d) - \sum_{j\in J_\infty} \ord_{\infty}R(Y_j),$$
so that
$$\deg ( Q_2(z) \det M(z)) \leq u + \deg Q_2 \leq \sum_{\sigma\in\Sigma\setminus \{\infty\} }  \sum_{j\in J_\sigma}  \ord_{\sigma}R(Y_j) + \tau + \cstun$$
using Eq. \eqref{eqhypdetnn}, where $\cstun$ depends only on $A$ and $\Sigma$ (since $0\leq \Card J_\sigma \leq q$ for any $\sigma$). To sum up,
 we have found a polynomial $Q_1$ of degree at most $\tau + \cstun$ such that
$$\det M(z) = \frac{Q_1(z)}{Q_2(z)} \prod _{\sigma\in\Sigma\setminus \{\infty\}} (z-\sigma)^{ \sum_{j\in J_\sigma}  \ord_{\sigma}R(Y_j)}.$$

\bigskip

\noindent {\bf Step 3: } Evaluation at $\alpha$.

\medskip

Let $q_\alpha = \Card J_\alpha$ and $q'_\alpha  = q-q_\alpha$, where $J_\alpha = \emptyset $ if $\alpha\not\in\Sigma$; for simplicity we assume that $J_\alpha = \{1,\ldots,q_\alpha\}$. Since the solutions $Y_1$, \ldots, $Y_{q_\alpha}$ of the differential system $Y'=AY$ are linearly independent over $\C$, there exist solutions $Y_{q_\alpha+1}$, \ldots, $Y_q$ such that $(Y_1,\ldots, Y_q)$ is a local basis of solutions at $\alpha$. Let $\calY \in M_q(\loga)$ be the matrix with columns $Y_1$, \ldots, $Y_q$, where $\loga =  \C[\log(z-\alpha)][[z-\alpha]]$. Then$\tra \calY M$ is the matrix 
$[R(Y_i)^{(k-1)}]_{1\leq i,k\leq q}$.

For any subset $\AAA$ of $\unq$ of cardinality $q'_\alpha = q-q_\alpha$, we denote by $\Delta_\AAA$ the determinant of the submatrix of $[R(Y_i)^{(k-1)}] $ obtained by considering only the rows with index $i\geq q_\alpha+1$ and the columns with index $k\in \AAA$, and by $\Deltati_\AAA$ the one 
obtained by removing these rows and columns. Then Laplace expansion by complementary minors yields
\begin{equation} \label{eqLaplace}
\det \calY (z) \cdot  \det M(z) = \sum_{\AAA \subset \unq \atop \Card \AAA = q'_\alpha} \eps_\AAA \Delta_\AAA(z)\Deltati_\AAA(z)
\eneq
with $\eps_\AAA\in\{-1,1\}$. Now $\det \calY $ is the wronskian of $Y_1$, \ldots, $Y_q$: it is a solution of the first order differential equation
\begin{equation} \label{eqwr}
w'(z) = w(z) {\rm trace}(A(z))   .
\eneq
Moreover it is non-zero, and belongs to $\loga$. Therefore we have $\det \calY(z) \sim \gammaun (z-\alpha)^\oun$ as $z\to\alpha$, for some $\gammaun\in\C\etoile$ and  $\oun\in\N$. On the other hand we have $\det M(z) \sim \gammade (z-\alpha)^\ode$ with  $\gammade\in\C\etoile$ and  $\ode\in\N$ using Step 1 and the assumption that all entries of $M(z)$ are holomorphic at $\alpha$. Now for $f\in \loga\setminus\{0\}$ let $\ord_\alpha f$ denote its generalized order at $\alpha$, namely the maximal integer $N$ such that $f \in (z-\alpha)^N \loga$. Then Eq. \eqref{eqLaplace} shows that, for some subset $\AAA$,
\begin{equation} \label{eqmajootr}
\otr := \ord_\alpha \Delta_\AAA(z) \leq \oun + \ode - \ord_\alpha \Deltati_\AAA(z) \mbox{ and } 
 \Delta_\AAA(z) \sim \gammatr (z-\alpha)^\otr \mbox { with }   \gammatr\in\C\etoile.
 \eneq
Now letting $\omega_\alpha =\sum_{j\in J_\alpha} \ord_\alpha R(Y_j)$  if $\alpha \in\Sigma$ and 
 $\omega_\alpha = 0$ otherwise, Step 2 shows that $\ode\leq \omega_\alpha + \tau + \cstun$. Moreover the order of vanishing at $\alpha$ of any non-zero solution of Eq. \eqref{eqwr}, and in particular $\oun$,  can be bounded from above in terms of $A$ only. At last, for any $i \in J_\alpha= \{1,\ldots,q_\alpha\}$ and any $k\in\unq$ the function $R(Y_i)^{(k-1)}$ vanishes at $\alpha$ with order at least $\ord_\alpha R(Y_i) - (q-1)$ so that 
 $\ord_\alpha \Deltati_\AAA(z) \geq \omega_\alpha - q_\alpha(q-1)$. Therefore Eq. \eqref{eqmajootr} yields
$\otr \leq  \tau + \csttr$ for some constant $\csttr$ depending only on $A$ and $\Sigma$. 

Now let us consider the $\otr$-th derivative $\Delta_\AAA^{(\otr)}(z)$: it has a finite non-zero limit (equal to $\gammatr \otr!$) as $z\to \alpha$. Moreover this derivative is a $\Z$-linear combination of determinants of matrices of the form
$$N_{k_1,\ldots,k_{q'_\alpha}} = [R(Y_{q_\alpha + i})^{(k_j-1)}]_{1\leq i,j \leq q'_\alpha}$$
with $1\leq k_1 < \ldots < k_{q'_\alpha} \leq q + \otr \leq \tau + \cstnum$; this constant $\cstnum$ 
(which depends only on $A$ and $\Sigma$)
 is the one in the statement of assertion $(ii)$ of Theorem \ref{thzeronum}. Now $P_{k,i}$ is assumed to be holomorphic at $\alpha$ for any $i$ and any $k\leq \tau + \cstnum$, so that $
R(Y_i)^{(k-1)} \in \loga$ since $Y_i \in M_{q,1}(\loga)$. Therefore $\det N_{k_1,\ldots,k_{q'_\alpha}} \in \loga$; since 
 $\Delta_A^{(\otr)}(z)$ has a finite  non-zero limit as $z\to \alpha$,  there exists at least one term $\det N_{k_1,\ldots,k_{q'_\alpha}} $ in the above-mentioned $\Z$-linear combination which also has a finite  non-zero limit as $z\to \alpha$. For this tuple we consider the equality $\tra \calYti \Mti = 
 N_{k_1,\ldots,k_{q'_\alpha}}$, where $\calYti  \in M_{q,q'_\alpha}(\loga)$ is the matrix with columns $Y_{q_\alpha+1}$, \ldots, $Y_q$, and $\Mti = [P_{k_j,i}]_{1\leq i \leq q, 1\leq j \leq q'_\alpha}$. The Cauchy-Binet formula  yields
\begin{equation} \label{eqCB1} 
\det N_{k_1,\ldots,k_{q'_\alpha}} = \sum_{B \subset \unq \atop \Card B = q'_\alpha} \det \tra \calYti_B \cdot \det \Mti_B
\eneq
where $\calYti_B$ (resp. $\Mti_B$) is the square matrix consisting in the rows of $\calYti$ (resp. of $\Mti$) corresponding to indices in $B$.

Let $\evalpha : \loga\to\C$ denote regularized evaluation at $\alpha$, defined by $\evalpha (f) = c_{0,0}$ for any $f = \sum_{i,j} c_{i,j} (z-\alpha)^i (\log(z-\alpha))^j$. The important point is that $\evalpha$ is a $\C$-algebra homomorphism, and that $\evalpha(f)$ is equal to the limit of $f(z)$ as $z\to\alpha$ whenever this limit exists. Extending $\evalpha$ coefficientwise to matrices, Eq. \eqref{eqCB1} yields
$$\evalpha \Big(\det N_{k_1,\ldots,k_{q'_\alpha}}\Big) = \sum_{B \subset \unq \atop \Card B = q'_\alpha} \evalpha \Big(\det \tra \calYti_B \Big)\cdot \evalpha \Big(\det \Mti_B\Big).$$
Now the left hand side is non-zero, so that $\evalpha (\det \Mti_B) \neq 0$ for some $B$. Since all coefficients $P_{k,i}$ are holomorphic at $\alpha$, so is $\det \Mti_B $ and therefore $\det (\Mti_B(\alpha) )=  \evalpha (\det \Mti_B) \neq 0$. We have found an invertible submatrix of $M(\alpha)$ of size $q'_\alpha$, so that $\rk M(\alpha) \geq q'_\alpha$: this concludes the proof of Theorem \ref{thzeronum}.

\section{Diophantine part  of the proof} \label{secdio}

In this section we prove Theorem \ref{thminocaractnv} stated in the introduction, and give  in details  new proofs of the Ball-Rivoal theorem and Nishimoto's Theorem \ref{thminocaract}. To provide  a unified treatment, we state a general result (namely Theorem \ref{thprinc}) and deduce these results from it in \S \ref{subsec31nv}. In order to help the reader, we first sketch the proof of Theorem \ref{thprinc} in \S \ref{subsec32}, then construct the linear forms (\S \ref{subsec33}), apply the zero estimate (namely Theorem \ref{thzeronum}) to obtain in invertible matrx (\S \ref{subsec34}), and study the arithmetic and asymptotic properties (\S \ref{subsec35}). At last we state and prove Siegel's linear independence criterion in \S \ref{subsec36}.

\subsection{Statement of the main theorem and consequences} \label{subsec31nv}

\begin{Th} \label{thprinc}
Let $N \geq 1$, and $f: \N \to \C$ be such that $f(n+N)=f(n)$ for any $n$. Let $p\in\{0,1\}$, $a\geq 2$, and $z_0 \in\{1, e^{i\pi/N}\}$; put
$$\xi_j = \sum_{n=1}^\infty \frac{f(n)z_0^n}{n^j} \mbox{ for any } j \in\una,$$
except that $\xi_1 = 0$ if $z_0=1$. Then as $a\to\infty$,
$$\dim_\Q\Span_\Q(\{\xi_j, \, \, 1\leq j\leq a, \, j\equiv p \bmod 2\} ) \geq \frac{1+o(1)}{N +\log 2}\log a.$$
\end{Th}

We refer to \S \ref{secBR} for the special case of the Ball-Rivoal theorem.

\bigskip

Let us deduce Theorems \ref{thminocaract} and  \ref{thminocaractnv} stated in the introduction from this result. Let $\chi$ be a Dirichlet character mod $d$. Its conductor is the smallest  divisor $e$  of $d$ for which there exists a character $\chi'$ mod $e$ such that $\chi(n) = \chi'(n)$ for any $n$ coprime to $d$. Comparing the $L$-functions of $\chi$ and $\chi'$  (see  for instance \cite[\S\S 3.2 and 3.3]{IwaniecKowalski}) yields
$$L(\chi,s)  = L(\chi',s)  \prod_{p | d \atop p \not\vert e} (1 - \chi'(p) p^{-s})$$
so that $\delta_{\chi,p,a} = \delta_{\chi ' ,p,a} $ for any $p$, $a$ (with the notation of  Theorem \ref{thminocaract}). Therefore we may assume that $e=d$, i.e. $\chi$ is primitive. Then Theorem \ref{thminocaract} follows from Theorem \ref{thprinc} by letting $z_0=1$ and $f=\chi$.

\bigskip

To prove Theorem \ref{thminocaractnv}, we first prove that for any primitive Dirichlet character $\chi$ modulo a multiple $e$ of 4, 
\begin{equation} \label{eqchisym}
\chi(n+\frac{e}2) = - \chi(n) \mbox{ for any } n\in\Z.
\eneq
Indeed we have $n(\frac{e}2+1) \equiv n+\frac{e}2 \bmod e$ if $n$ is odd, so that $\chi(n+\frac{e}2) =  \chi(n) \chi(\frac{e}2+1)$ for any $n\in\Z$ (since both sides vanish if $n$ is even). Moreover $(\chi(\frac{e}2+1))^2=1$ since $ (\frac{e}2+1)^2  \equiv 1 \bmod e$, and  $ \chi(\frac{e}2+1) \neq 1$ because $\chi$ is primitive (so that $\chi(n+\frac{e}2) \neq \chi(n) $ for some $n$). Therefore  $ \chi(\frac{e}2+1) =-1$: this concludes the proof of \eqref{eqchisym}. 

Now let $N = e/2$ and define $f: \N \to\C $ by $f(r)  = \chi(r) z_0^{-r}$ for any $r\in\{1,\ldots,N\}$, where $z_0 = e^{i\pi/N}$. Then Eq. \eqref{eqchisym} yields
$$\sum_{n=1}^\infty  \frac{f(n)z_0^n}{n^j} = \sum_{r=1}^N f(r) z_0^r  \sum_{n\geq 1\atop n\equiv r \bmod 2N}\Big( \frac{1}{n^j}-\frac1{(n+N)^j}\Big) = 
\sum_{n=1}^\infty  \frac{\chi(n)}{n^j} = L(\chi,j)$$
so that  Theorem \ref{thprinc} implies  Theorem \ref{thminocaractnv}.

\subsection{Sketch of the proof} \label{subsec32}

To prove Theorem \ref{thprinc}, we let $r,n\geq 1$ be such that $r < \frac{a}{2N}$ and $N$ divides $n$. We define $\xi'_1$, \ldots, $\xi'_{a+N}$
as follows:
\begin{equation} \label{eqdefxip}
\left\{
\begin{array}{cll}
\xi'_j &=  \hspace{0.3cm} 2(-1)^p \xi_j &\mbox{  for $j \in \una $ such that } j \equiv p \bmod 2 \\
\xi'_j &=  \hspace{0.3cm}  0 &\mbox{  for $j \in \una $ such that  } j \not\equiv p \bmod 2 \\
\xi'_{a+1+\lambda} &=   \hspace{0.3cm} z_0^\lambda f(\lambda)  &\mbox{  for any } \lambda\in\{0,\ldots,N-1\}.
\end{array}
\right.
\eneq
We also let 
$$\delta_n = (N d_n)^a N^{an/N} , 
$$
and define $\ttt$ to be equal to 1 if $z_0 = e^{i\pi/N}$, and equal to 2 otherwise (i.e., if $z_0=1$). 

In \S \ref{subsec33}  (see \eqref{eqdefski}) we shall construct integers $s_{k,i}$, $\ttt \leq i \leq a+N$, such that as $n\to+\infty$:
\begin{equation} \label{eqestilamk}
\max_{\ttt \leq i \leq a+N} | s_{k,i} | \leq \beta^{n(1+o(1))} \mbox { and }
\Big|  \sum_{i=\ttt} ^{a+N}  s_{k,i} \xi'_i \Big|   \leq \alpha ^{n(1+o(1))}
\eneq
 where 
$$ \alpha = e^a 4^{a/N - r} (N+1)^{2r+2} r^{-a/N + 4r+2} \mbox { and }
 \beta = (2e^N)^{a/N} (rN+1)^{2r+2}.$$
Then Lemma \ref{lemmatinv} (that will be stated and proved in \S \ref{subsec34} using Theorem \ref{thzeronum}) provides a positive constant $\cstnum$ (which depends only on $a$ and $N$) and integers $1 \leq k_{\ttt} < k_{\ttt+1} < \ldots < k_{a+N}\leq \cstnum$ (which depend on $a$, $N$, $r$, and $n$) such that the matrix $[s_{k_j,i}]_{\ttt\leq i,j\leq a+N}$ is invertible.  Since $k_j \leq \cstnum$ for any $j$, the symbols $o(1)$ in \eqref{eqestilamk} with $k=k_j$ can be made uniform with respect to $k$. Therefore Siegel's linear independence criterion applies
 (see \S \ref{subsec36}). Taking  $a$ very large, $N$ fixed, and $r$ equal to the integer part of $\frac{a}{(\log(a))^2}$  concludes  the proof of Theorem \ref{thprinc} since $\xi'_i =0$ if $i\leq a$ and $i\not\equiv p \bmod 2$, and 
$$1 - \frac{\log\alpha}{\log\beta} = \frac{1+\eps_a}{N+\log 2} \log a \mbox{ where } \lim_{a\to+\infty} \eps_a = 0.$$

\subsection{Construction of the linear forms}\label{subsec33}

Let $a$, $r$, $N$ be positive integers such that $1 \leq r < \frac{a}{2N}$. For any integer multiple $n$ of $N$ we let
$$F(t) = (n/N)!^{a-2rN} \frac{(t-rn)_{rn} (t+n+1)_{rn}}{\prod_{h=0}^{n/N}(t+Nh)^a}.$$
Then $F$ is a rational function, and its degree $-\dz$ satisfies 
\begin{equation} \label{eqdefdz}
\dz := a(\frac{n}{N} +1) - 2rn = -\deg F \geq n+a\geq 2.
\eneq
Its partial fraction expansion reads
$$F(t) = \sum_{h=0}^{n/N}\sum_{j=1}^a \frac{p_{j,h}}{(t+Nh)^j}$$
with rational coefficients $p_{j,h}$. Let
$$P_j(z) = \sum_{h=0}^{n/N} p_{j,h} z^{Nh} \in\Q[z]_{\leq n} \mbox{ for any } j\in\{1,\ldots,a\},$$
and also
\begin{equation} \label{eqdefszinf}
S_0(z) = \sum_{t=n+1}^\infty F(-t) z^t, \hspace{2cm} 
S_\infty(z) = \sum_{t=1}^\infty F( t) z^{-t} .
\eneq
As in \cite{BR} we have
$$S_\infty(z) = V(z) + \sum_{j=1}^a P_j(z) \Li_j(1/z)$$
where 
$$V(z) =  -\sum_{t=0}^{n-1} z^t \sum_{j=1}^a \sum_{h= \lceil (t+1)/N\rceil} ^{n/N}  \frac{p_{j,h}}{(Nh-t)^j} \in\Q[z]_{\leq n}.$$
In the same way (see \cite{FR}) we have 
$$S_0(z) = U(z) + \sum_{j=1}^a P_j(z) (-1)^j \Li_j(z)$$
with the same polynomials $P_1$, \ldots, $P_a$, and 
$$U(z) =  -\sum_{t=1}^n z^t \sum_{j=1}^a \sum_{h=0}^{\lfloor (t-1)/N\rfloor} \frac{p_{j,h}}{(t-Nh)^j} \in\Q[z]_{\leq n}.$$

\bigskip

Now let $P_{1,j} = P_j$ for any $j\in\una$, and define inductively $P_{k,j}\in\Q(z)$ by
\begin{equation} \label{eqdefpkj}
P_{k,j}(z) = P'_{k-1,j}(z) - \frac1{z} P_{k-1,j+1}(z) \mbox{ for any } k\geq 2 \mbox{ and any } j\in\una,
\eneq
where   $P_{k-1,a+1} = 0$ for any $k$; we shall check in \S \ref{subsec34} below that this notation $P_{k,j}$ is consistent with the one used in the introduction. We let also $U_1 = U$, $V_1 = V$, and define $U_k$, $V_k$ for any $k\geq 2$ by the recurrence relations
\begin{equation} \label{eqdefpz}
U_k(z) = U'_{k-1}(z) - \frac1{1-z} P_{k-1,1}(z),
\eneq
\begin{equation} \label{eqdefpinf}
V_k(z) = V'_{k-1}(z) + \frac1{z(1-z)} P_{k-1,1}(z).
\eneq
Then for any $k\geq 1$ we have
\begin{equation} \label{eqszeroder}
S_0^{(k-1)}(z) = U_k(z) + \sum_{j=1}^a P_{k,j}(z) (-1)^j \Li_j( z)
\eneq
\begin{equation} \label{eqsinfder}
\mbox{ and } S_\infty^{(k-1)}(z) = V_k(z) + \sum_{j=1}^a P_{k,j}(z)  \Li_j(1/z).
\eneq

\bigskip

Moreover Eqns. \eqref{eqdefpkj}, \eqref{eqdefpz} and \eqref{eqdefpinf} show that the rational functions $P_{k,j}$ with $1\leq j \leq a$ (resp. $U_k$ and $V_k$) have only 0 (resp. only 0 and 1) as possible finite poles. Now we have 
$$S_0^{(k-1)}(z) = \sum_{t=n+1} ^\infty F(-t)(t-k+2)_{k-1} z^{t-k+1} \mbox{ for } |z| < 1 $$
$$\mbox{ and } S_\infty^{(k-1)}(z) =  \sum_{t=1} ^\infty F( t) (-1)^{k-1} (t)_{k-1} z^{-t-k+1} \mbox{ for } |z| > 1.$$
Let us assume that $k-1\leq \dz-2$, where $\dz = -\deg F$ is defined by Eq. \eqref{eqdefdz};  then these formulas hold also when $|z|=1$ and  we may let $z$ tend to 1 in Eqns. \eqref{eqszeroder} and \eqref{eqsinfder}. Since $P_{k,j} $ is holomorphic at $z=1$ for any $k\geq 1$ and any $j$, a possible   divergence may come only from  poles of $U_k$ or $V_k$ at $z=1$, or from the logarithmic term involving $\Li_1(z)$ or $\Li_1(1/z)$. Since a pole and a  logarithmic term cannot cancel each other out, and $S_0^{(k-1)}(z) $ and $S_\infty^{(k-1)}(z) $ have finite limits as $z\to 1$, we obtain:
\begin{equation} \label{eqokenun}
\mbox{For any }k\leq \dz-1,  \hspace{.7cm} P_{k,1}(1) = 0  \hspace{.2cm} \mbox{ and }  \hspace{.2cm} U_k, V_k \mbox{ do not have a pole at } z=1.
\eneq

\bigskip

Now let $k\leq \dz-1$, and $z\in\C$ be such that $|z|=1$. Then Eqns. \eqref{eqszeroder} and  \eqref{eqsinfder} hold, upon agreeing that the sums start at $j=2$ if $z=1$; the same remark applies in what follows. Since $P_j(z)\in\Q[z^N]$ for any $j\in\una$, Eq. \eqref{eqdefpkj} yields $P_{k,j} \in z^{1-k} \Q[z^N]$ (see the proof of Proposition \ref{propformelin} in \S \ref{subsec35} for details). On the other hand, since 
 $U_k,V_k \in \Q[z,z^{-1}]$ for any $k\leq \dz-1$, we can write
\begin{equation} \label{eqpkzinfom}
z^{ k-1 } U_k(z) = \sum_{\lam = 0}^{N-1} z^\lam U_{k,\lam}(z) 
\mbox{ and }
z^{ k-1 }V_k(z) = \sum_{\lam = 0}^{N-1} z^\lam V_{k,\lam}(z) 
\eneq
with $U_{k,\lam},  V_{k,\lam} \in  \Q[z^N,z^{-N}]$.  Then  Eqns. \eqref{eqszeroder} and  \eqref{eqsinfder} yield
\begin{equation} \label{eqszeroderbis}
z^{k-1} S_0^{(k-1)}(z) =\sum_{\lam = 0}^{N-1} z^\lam U_{k,\lam}(z)  + \sum_{j=1}^a z^{k-1}  P_{k,j}(z) (-1)^j \Li_j( z)
\eneq
\begin{equation} \label{eqsinfderbis}
\mbox{ and } z^{k-1}  S_\infty^{(k-1)}(z) = \sum_{\lam = 0}^{N-1} z^\lam V_{k,\lam}(z)  + \sum_{j=1}^a z^{k-1}  P_{k,j}(z)  \Li_j(1/z).
\eneq
The point now is that $U_{k,\lam}(z)  $, $V_{k,\lam}(z)  $, and $z^{k-1}  P_{k,j}(z)$ depend only on $z^N$. For any $\ell\in\unN$ we consider 
\begin{equation} \label{eqdefmunu}
\mu_\ell = \frac1{N} \sum_{\lam = 1}^N f(\lam) \om^{-\ell \lam}.
\eneq
Let $z_0 \in \{1, e^{i\pi/N}\}$ and $p\in\{0,1\}$ be as in Theorem \ref{thprinc}, and recall that $\om = e^{2i\pi/N}$. For any $k\leq \dz-1$ we let
$$\Lambda_k = \sum_{\ell=1}^N \mu_\ell \Big[ (\om^\ell z_0)^{k-1} S_0^{(k-1)}(\om^\ell z_0) + (-1)^p   (\om^\ell z_0)^{1-k} S_\infty^{(k-1)}(\frac1{\om^\ell z_0})\Big].$$
Then Eqns.  \eqref{eqszeroderbis} and  \eqref{eqsinfderbis} yield,  since $U_{k,\lam}(z)  $, $V_{k,\lam}(z)  $, and $z^{k-1}  P_{k,j}(z)$ depend only on $z^N$ and $(\om^\ell z_0)^{ N} =(\om^\ell z_0)^{-N} = z_0^N$:
\begin{eqnarray*}
\Lambda_k &=& \sum_{\lam = 0}^{N-1}  \Big[ \Big(  \sum_{\ell=1}^N \mu_\ell (\om^\ell z_0)^{\lambda}\Big) U_{k,\lam}(z_0) + (-1)^p    \Big(  \sum_{\ell=1}^N \mu_\ell (\om^\ell z_0)^{-\lambda}\Big) V_{k,\lam}(z_0) \Big] \\
&&+ \sum_{j=1}^a z_0^{k-1}  P_{k,j}(z_0)   \sum_{\ell=1}^N \mu_\ell 
 \Li_j( \om^\ell z_0) ( (-1)^j + (-1)^p).
\end{eqnarray*}
  Now Eq. \eqref{eqdefmunu} yields
$$  \sum_{\ell=1}^N \mu_\ell \om^{n  \ell} = f(n) \mbox{ for any $n\in\Z$, so that }
  \sum_{\ell=1}^N \mu_\ell \Li_j(\om^\ell z_0) = \sum_{n=1}^\infty \frac{f(n) z_0^n}{n^j} =  \xi_j \mbox{ for any } j \leq a.$$ 
Letting $V_{k,N} = V_{k,0}$ we obtain: 
$$\Lambda_k = 2(-1)^p \sum_{1\leq j \leq a \atop j\equiv p \bmod 2}z_0^{k-1} P_{k,j}(z_0) \xi_j + \sum_{\lam=0}^{N-1} (U_{k,\lam}(z_0)+(-1)^p V_{k,N-\lam}(z_0))z_0^\lam f(\lam).$$
As announced in \S \ref{subsec32} we now define the coefficients $s_{k,i}$:
\begin{equation} \label{eqdefski}
\left\{
\begin{array}{l}
s_{k,i} =  \delta_n z_0^{k-1}P_{k,i}(z_0) \mbox{ for } 1\leq i \leq a, \\
s_{k,a+1+\lam} = \delta_n  (U_{k,\lam}(z_0)+(-1)^p V_{k,N-\lam}(z_0))  \mbox{ for } 0\leq \lam \leq N-1,
\end{array}
\right.
\eneq
where
$\delta_n = (N d_n)^a N^{an/N} $, 
so that 
$$\delta_n \Lambda_k=  \sum_{i = \ttt } ^{a+N}   s_{k,i}  \xi'_i $$
since $\xi'_1 = 0$ if $z_0 = 1$ (recall from \S \ref{subsec32} that $\ttt=2$ in this case, and $\ttt=1$ otherwise, i.e. if $z_0 = e^{i\pi/N}$; $\xi'_i$ is defined in Eq. \eqref{eqdefxip}). 

Since $z_0^N\in\{-1,1\}$ and $z^{k-1}P_{k,j}(z)$, $U_{k,\lam}(z)$ and $V_{k,N-\lam}(z)$ are polynomials in $z^N$ with rational coefficients, the  numbers $s_{k,1}$, \ldots, $s_{k, a+N}$ are rational. We shall prove in Proposition \ref{propformelin} (\S \ref{subsec35})   that they are integers,  thanks to the factor $\delta_n$.

\subsection{Application of the zero estimate}\label{subsec34}

\newcommand{\Pbar}{\overline P}

\newcommand{\Pgras}{{\bf P}}

In this section we  deduce from Theorem \ref{thzeronum} the following lemma, used at the end of \S \ref{subsec32}. It provides an invertible matrix which enables us to apply Siegel's linear independance criterion (see \S \ref{subsec36}). 

\begin{Lem} \label{lemmatinv} In the setting of \S \ref{subsec32}, let $\ttt = 1$ if $z_0  = e^{i\pi/N}$ and $\ttt=2$ if $z_0=1$; let $s_{k,i}$ be defined by Eq. \eqref{eqdefski}. Then there exist  a positive constant $\cstnum$ (which depends only on $a$ and $N$) and integers $1 \leq k_{\ttt} < k_{\ttt+1} < \ldots < k_{a+N}\leq \cstnum$ (which depend on $a$, $N$, $r$, and $n$) such that the matrix $[s_{k_j,i}]_{\ttt\leq i,j\leq a+N}$ is invertible. 
\end{Lem}

To begin with, let us recall from \S  \ref{subsec32}  that $\om = e^{2i\pi/N}$, 
  $a$, $r$, $N$, $n$ are positive integers such that $1 \leq r < \frac{a}{2N}$,    $n$ is a   multiple  of $N$, and 
  $$F(t) = (n/N)!^{a-2rN} \frac{(t-rn)_{rn} (t+n+1)_{rn}}{\prod_{h=0}^{n/N}(t+Nh)^a}.$$
We have
$$S_0(z) = \sum_{t=n+1}^\infty F(-t) z^t =  U(z) + \sum_{j=1}^a P_j(z) (-1)^j \Li_j(z)$$
$$\mbox{ and } S_\infty(z) = \sum_{t=1}^\infty F( t) z^{-t} =  V(z) + \sum_{j=1}^a P_j(z) \Li_j(1/z).$$
Since $P_j \in \C[z^N]$ for any $j\in\una$, we have $P_j(\om^\ell z) = P_j(z)$ for any $\ell\in\Z$. Therefore letting
\begin{equation} \label{eqc111}
R_{0,\ell}(z) = S_0(\om^\ell z),
\hspace{0.6cm} 
R_{\infty,\ell}(z) = S_\infty(\om^\ell z),
\hspace{.6cm} 
\Pbar_{0,\ell}(z) = U(\om^\ell z),
\hspace{.6cm} 
\Pbar_{\infty,\ell}(z) = V(\om^\ell z) 
\eneq
for any $\ell\in\unN$, we have
\begin{equation} \label{eqP1}
R_{0,\ell}(z) = \Pbar_{0,\ell}(z) +   \sum_{j=1}^a P_j(z) (-1)^j \Li_j(\om^\ell z) = O(z^{(r+1)n+1}), \hspace{0.5cm} z\to 0,
\eneq
\begin{equation} \label{eqP2}
\mbox{ and }
R_{\infty,\ell}(z) = \Pbar_{\infty,\ell}(z) +   \sum_{j=1}^a P_j(z)   \Li_j(\frac1{\om^\ell z}) = O(z^{ -rn-1}), \hspace{0.5cm} z\to \infty.
\eneq
Moreover, recall that $\dz = -\deg F = a(\frac{n}{N}+1)-2rn$;  Lemma 3 of \cite{FR} shows that
$$
  \sum_{j=1}^a P_j(z)  (-1)^{j-1}\frac{(\log z)^{j-1}}{(j-1)!} = O((z-1)^{\dz-1}),  \hspace{0.5cm} z\to 1.$$
Using again the fact that $  P_j(\om^{-\ell}  z) = P_j(z)$, we obtain for any $\ell \in\unN$:
\begin{equation} \label{eqP3}
R_{\om^\ell  }(z) :=    \sum_{j=1}^a P_j(z)  (-1)^{j-1}\frac{(\log (\om^{-\ell}z)^{j-1}}{(j-1)!} = O((z-\om^\ell)^{\dz-1}),  \hspace{0.5cm} z\to \om^\ell.
\eneq
Combining Eqns. \eqref{eqP1},  \eqref{eqP2},  and  \eqref{eqP3} with $1\leq \ell \leq N$, we have solved a simultaneous Pad\'e approximation problem. The $(n+1)(a+2N)$ unknowns are the coefficients of $P_1$, \ldots, $P_a$,  $\Pbar_{0,1}$, \ldots, $\Pbar_{0,N}$, $\Pbar_{\infty,1}$, \ldots, $\Pbar_{\infty,N}$, which are polynomials of degree less than or equal to $n$. There are
$$2N((r+1)n+1)+N(\dz-1) = n(a+2N)+(a+1)N$$
linear equations, since a priori we have $R_{\infty, \ell}(z) = O(z^n)$ as $z\to\infty$. The difference between the number of unknowns and the number of equations is equal to $N-a(N-1)$. If $N=1$ this is equal to 1: the Pad\'e approximation problem is exactly  \eqref{eqpadeFR}, i.e. the one of  \cite[Th\'eor\`eme 1]{FR}, which has a unique solution up to proportionality. Whenever $N\geq 2$ we have $N-a(N-1) < 0$: the problem we have solved has more equations than unknowns. This is due to the fact that we always assume $n$ to be an integer multiple of $N$. Anyway to complete the proof, it is sufficient to bound from above the difference between the number of unknowns and the number of equations by a constant independent from $n$; we do not need to study whether the Pad\'e  approximation problem has a unique solution or not. 

\bigskip

Let $q = a+2N$, and $A\in M_q(\C(z))$ be the matrix of which the coefficients $A_{i,j}$ are given by:
$$
\left\{
\begin{array}{l}
A_{i,i-1}(z) = \frac{-1}{z} \mbox{ for any } i \in \{2,\ldots,a\}\\
A_{1,a+ \ell}(z) = \frac{ \om^\ell}{ \om^\ell z-1 } \mbox{ for any } \ell\in\unN \\
A_{1,a+N+\ell}(z) = \frac{1}{z(1-\om^\ell z  ) } \mbox{ for any } \ell\in\unN
\end{array}
\right.
$$
and all other coefficients are zero. We consider the following solutions of the differential system $Y'=AY$, with $1\leq\ell\leq N$:
$$Y_{0,\ell}(z) = \tra \Big( -\Li_1(\om^\ell z), \Li_2(\om^\ell z), \ldots, (-1)^a \Li_a(\om^\ell z), 0,  \ldots, 0,1,0, \ldots,0\Big),$$
$$Y_{\infty,\ell}(z) = \tra \Big(  \Li_1(\frac1{\om^\ell z}), \Li_2(\frac1{\om^\ell z}), \ldots,   \Li_a(\frac1{\om^\ell z}), 0,  \ldots, 0,1,0, \ldots,0\Big),$$
$$Y_{\om^\ell}(z) = \tra \Big( 1, -\log(\om^{-\ell}z), \frac{(\log(\om^{-\ell}z))^2}{2!}, \ldots,  (-1)^{a-1} \frac{(\log(\om^{-\ell}z))^{a-1}}{(a-1)!},0 , \ldots,0\Big)$$
where the coefficient 1 in $Y_{0,\ell}(z) $ (resp. $Y_{\infty,\ell}(z) $) is in position $a+\ell$ (resp. $a+N+\ell$). 

We let $J_0 = \{(0,1),(0,2),\ldots,(0,N)\}$, $J_\infty = \{(\infty,1),(\infty,2),\ldots,(\infty,N)\}$, $J_{\om^\ell}  = \{\om^\ell\}$ for $1\leq \ell \leq N$, and  $\Sigma = \{0,\infty\} \cup \{\om^\ell, 1\leq \ell \leq N\}$. We also let $P_{a+\ell}(z) = \Pbar_{0,\ell}(z) = U(\om^\ell z)$ and $P_{a+N+\ell} = \Pbar_{\infty,\ell}(z) = V(\om^\ell z) $ for any $\ell\in\unN$. Then with the notation of the introduction we have $R(Y_{0,\ell}) = R_{0,\ell}(z)$, $R(Y_{\infty,\ell}) = R_{\infty,\ell}(z)$, and $R(Y_{\om^\ell}) = R_{\om^\ell}(z)$  for any $\ell\in\unN$.

Since $P_a$ is not the zero polynomial, we have $R_{\om^\ell}(z)\neq 0$ for any $\ell$; the $\C$-linear independence of $R_{0,1}(z)$, \ldots, $R_{0,N}(z)$  (resp. of $R_{\infty,1}(z)$, \ldots, $R_{\infty,N}(z)$) follows directly  (resp. up to changing $z$ to $1/z$) from the following lemma, which is not difficult to prove using monodromy (see \cite{MonodromiePolylogs}).

\begin{Lem} 
The functions 1 and $\Li_j(\om^\ell z)$, for $j\geq 1$ and $1\leq \ell \leq N$, are linearly independent over $\C(z)$.
\end{Lem}

Eqns. \eqref{eqP1}, \eqref{eqP2}, and \eqref{eqP3} yield $\ord_0 (R_{0,\ell}(z)) \geq (r+1)n+1$,  $\ord_\infty (R_{\infty,\ell}(z)) \geq rn+1$,  and 
 $\ord_{\om^\ell}  (R_{\om^\ell}(z)) \geq \dz-1$ for any $\ell\in\unN$, so that 
 $$\sum_{\sigma\in\Sigma} \sum_{j\in J_\sigma} \ord_\sigma R_j (z) \geq (2r+1)Nn+N(\dz+1) = (n+1) q-nN -\tau \mbox{ with } \tau =N- a(N-1); $$
 here $q = a+2N$, and we recall that $\dz = -\deg F = a(\frac{n}{N}+1)-2rn$. This number $\tau$ is exactly the 
 difference between the number of unknowns and the number of equations
computed after Eq. \eqref{eqP3}. 

 Now for any $k\geq 1$ and any $\ell\in\unN$ we let 
 \begin{equation} \label{eqdefpbar}
 \Pbar_{k,0,\ell} = \om^{\ell(k-1)}U_k(\om^\ell z) \mbox{ and } 
 \Pbar_{k,\infty,\ell} = \om^{\ell(k-1)}V_k(\om^\ell z),
 \eneq
 and
 $$\Pgras_k = \tra\Big( P_{k,1},P_{k,2}, \ldots,P_{k,a},  \Pbar_{k,0,1}, \ldots,  \Pbar_{k,0,N},  \Pbar_{k,\infty,1}, \ldots,  \Pbar_{k,\infty,N}\Big)\in M_{q,1}(\C(z)),$$
so that $\Pgras_1 = \tra (P_1,\ldots,P_{a+2N})$.
 Then it is not difficult to check that  
 $$\Pgras_k = \Big( \frac{\dd}{\dd z}  + \tra A \Big)^{k-1} \Pgras_1.$$
 To illustrate this equality, we notice that Eq. \eqref{eqc111} yields 
 $$R_{0,\ell}^{(k-1)} = \Pbar_{k,0,\ell}(z) +\sum_{j=1}^a P_{k,j} (z) (-1)^j \Li_j(\om^\ell z)$$
$$\mbox{  and }
  R_{\infty,\ell}^{(k-1)} = \Pbar_{k,\infty,\ell}(z) +\sum_{j=1}^a P_{k,j} (z)   \Li_j(\frac1{\om^\ell z}) $$
since (as in \cite[Chapter 3, \S 4]{Shidlovski})
 $$S_{0 }^{(k-1)} = U_k(z) +\sum_{j=1}^a P_{k,j} (z) (-1)^j \Li_j(  z)$$
$$\mbox{  and }
 S_{\infty }^{(k-1)} = V_k(z) +\sum_{j=1}^a P_{k,j} (z)   \Li_j(1/z ). $$
 
   \bigskip
   
Provided $n$ is large enough, we have checked all assumptions of Theorem \ref{thzeronum} (using, among others, Eq. \eqref{eqokenun}). We apply this result with $\alpha = z_0$; recall that $z_0 \in \{1,e^{i\pi/N}\}$. In the case $z_0=1$, we obtain positive integers  $k_2< \ldots < k_q \leq  \cstnum$ such that the matrix with columns $\Pgras_{k_2}(1)$, \ldots,   $\Pgras_{k_q}(1)$ has rank $q-1$. Now $P_{k,1}(1) = 0$ for any $k\leq  \cstnum$ (using   Eq. \eqref{eqokenun} since $n$ is large enough) so that the first row of this matrix is identically zero. Removing this row yields the following invertible matrix (with $z_0=1$ and $\ttt=2$):
\begin{equation} \label{eqmatrigr2}
\left[ \begin{array}{c}
[z_0^{k_j-1}  P_{k_j,i}(z_0)]_{\ttt\leq i \leq a, \ttt \leq j\leq q }\\
\hline
[z_0^{k_j-1}  \Pbar_{k_j,0,i}(z_0)]_{1\leq i\leq N, \ttt\leq  j\leq q}\\
\hline
[z_0^{k_j-1}  \Pbar_{k_j,\infty,i}(z_0)]_{1\leq i\leq N, \ttt\leq  j\leq q}
\end{array}\right] .
\eneq
If $z_0 = e^{i\pi/N} \not\in\Sigma $ then  Theorem \ref{thzeronum} provides directly  $k_1< \ldots < k_q \leq  \cstnum$  such that the matrix \eqref{eqmatrigr2} with $\ttt=1$ is invertible.

Now Eq. \eqref{eqpkzinfom} with $z=\om^\ell z_0$ yields, since $U_{k, \lambda}\in \Q[z^N, z^{-N}]$:
$$\om^{(k-1)\ell} z_0^{k-1} U_k(\om^\ell z_0 ) = \sum_{\lambda=0}^{N-1} \om^{\ell \lambda} z_0^\lambda U_{k, \lambda}(z_0)
\mbox{ for any }\ell\in\unN.$$
Therefore we have for any $\lambda\in\zeroNmu$:
\begin{equation} \label{equklamzz}
U_{k ,\lambda}(z_0) = \frac{z_0^{k-1-\lambda}}{N} \sum_{\ell=1}^N \om^{(k-1-\lambda)\ell} U_k(\om^\ell z_0 ) = 
\frac{z_0^{k-1-\lambda}}{N}\sum_{\ell=1}^N \om^{ -\lambda \ell} \Pbar_{k,0,\ell}(z_0)
\eneq
using Eq. \eqref{eqdefpbar}. Moreover the same relation holds with $V_{k ,\lambda}$ and $\Pbar_{k,\infty,\ell}$ for  $\lambda\in\zeroNmu$. We recall that $s_{k,i} $ was defined in Eq. \eqref{eqdefski}  (\S \ref{subsec33}) by 
\begin{eqnarray*}
&s_{k,i} =& \delta_n z_0^{k-1} P_{k,i}(z_0)\mbox{ for } 1\leq i \leq a, \\
\mbox{  and } &s_{k, a+1+\lambda} =& \delta_n( U_{k,\lambda}(z_0) + (-1)^p V_{k,N-\lambda}(z_0)) \mbox{ for } 0\leq \lambda\leq N-1.
\end{eqnarray*}
For any $\lambda\in\zeroNmu$ we deduce that
$$
s_{k, a+1+\lambda} = \frac{\delta_n}{N} z_0^{-\lambda}  \sum_{\ell=1}^N \om^{-\lambda\ell} z_0^{k -1}  \Pbar_{k ,0,\ell }(z_0) \pm
 (-1)^p \frac{\delta_n}{N} z_0^{-N+\lambda}  \sum_{\ell=1}^N \om^{ \lambda\ell} 
z_0^{k -1}  \Pbar_{k ,\infty,\ell }(z_0)$$
where $\pm$ is $+$ if $1\leq \lambda \leq N-1$, and $z_0^N$ if $ \lambda=0$; indeed $V_{k,N}=V_{k,0}$ satisfies the equation analogous to  Eq. \eqref{equklamzz} with $\lambda=0$, but not with $\lambda=N$ if $z_0 = e^{i\pi/N}$. 

Let $M = [m_{i,j}]_{\ttt \leq i \leq a+N, \ttt \leq j \leq a+2N}$ be the matrix defined by:
\begin{eqnarray*}
m_{i,i} = \delta_n \mbox{ for any } i\in\{\ttt, \ldots, a\},\\
m_{a+1+\lambda, a+\ell} = \frac{\delta_n}{N} z_0^{-\lambda} \om^{-\lambda\ell} \mbox{ for any } \lambda\in\zeroNmu \mbox{ and any } \ell\in\unN,\\
m_{a+1+\lambda, a+N+\ell} = \pm (-1)^p \frac{\delta_n}{N} z_0^{-N+\lambda} \om^{ \lambda\ell} \mbox{ for any } \lambda\in\zeroNmu \mbox{ and any } \ell\in\unN,
\end{eqnarray*}
and all other coefficients are zero. Then $M$ has rank $a+N+1-\ttt$; denoting by $P \in {\rm GL}_{a+2N+1-\ttt}(\C)$ the matrix \eqref{eqmatrigr2}, the matrix $MP$ has rank $a+N+1-\ttt$. Now $MP$ is exactly the matrix $[s_{k_j ,i}]_{\ttt\leq i,j\leq a+N}$. This concludes the proof of Lemma \ref{lemmatinv}.

\subsection{Arithmetic and Asymptotic Properties}\label{subsec35}

In this section we prove the following result,  used  in the proof of Theorem \ref{thprinc}; see \S \ref{subsec32} for the notation.

\begin{Prop} \label{propformelin}
Let 
\begin{equation} \label{eqdefalphabeta}
\alpha = e^a 4^{a/N - r} (N+1)^{2r+2} r^{-a/N + 4r+2} \mbox{ and }
\beta = (2e^N)^{a/N} (rN+1)^{2r+2}.
\eneq
Then we have $s_{k,i}\in\Z$ for any $i\in\{1,\ldots,a+N\}$ and any $k\leq \dz-1$, and as $n\to\infty$:
$$\Big| \sum_{i = \ttt} ^{a+N} s_{k,i} \xi'_i  \Big| \leq  \alpha ^{n(1+o(1))}, \hspace{1cm} 
\max_{1\leq i \leq a+N}  | s_{k,i} | \leq \beta  ^{n(1+o(1))}.$$
\end{Prop}

In this proposition and throughout this section, we denote by $o(1)$ any sequence that tends to 0 as $n\to\infty$; it usually depends  also on $a$, $r$, $N$, and $k$. 
When Proposition \ref{propformelin} is applied in the proof of Theorem \ref{thprinc} (see \S \ref{subsec32}), this dependence is not a problem since $a$, $r$, $N$ are fixed parameters and $k$ is bounded from above by $\cstnum$.
At last we  recall  that $d_{n}$ is the least common multiple of 1, 2, \ldots, $n$, and that 
$$\delta_n = (N d_n)^a N^{an/N}.$$

\bigskip

Let us start with a lemma, in which (as in \S \ref{subsec33})
$$F(t) = (n/N)!^{a-2rN} \frac{(t-rn)_{rn} (t+n+1)_{rn}}{\prod_{h=0}^{n/N}(t+Nh)^a}
= \sum_{h=0}^{n/N}\sum_{j=1}^a \frac{p_{j,h}}{(t+Nh)^j}.
$$

\begin{Lem} \label{lempjh}
For any $j\in\una$ and any $h \in \zeronsN$ we have
\begin{equation} \label{eqdenompjh}
(Nd_{n/N})^{a-j} N^{an/N} p_{j,h} \in \Z
\eneq
\begin{equation} \label{eqmajopjh}
\mbox{ and }
| p_{j,h} | \leq \Big( (2/N)^{a/N} (rN+1)^{2r+2}\Big)^{n(1+o(1))}
\eneq
where $o(1)$ is a sequence that tends to 0 as $n\to\infty$ and may depend also on $N$, $a$, and $r$.
\end{Lem}

\Dem of Lemma \ref{lempjh}:   We follow the approach of \cite{Habsieger} and \cite{Colmez} by letting
\begin{eqnarray*}
F_0(t)  &= \frac{(n/N)!}{\prod_{h=0}^{n/N} (t+Nh) } & = \sum_{h=0}^{n/N} \frac{(-1)^h N^{-n/N} \combitiny{n/N}{h}}{t+Nh} ,\\ 
G_i(t)  &= \frac{(t-in/N)_{n/N}}{\prod_{h=0}^{n/N} (t+Nh) } & = \sum_{h=0}^{n/N} \frac{(-1)^{h+n/N} N^{-n/N} \combitiny{n/N}{h}  \combitiny{Nh+in/N}{n/N}}{t+Nh}  \mbox{ for } 1\leq i \leq rN, \\ 
H_i(t)  &= \frac{(t+1+in/N)_{n/N}}{\prod_{h=0}^{n/N} (t+Nh) } & = \sum_{h=0}^{n/N} \frac{(-1)^{h } N^{-n/N} \combitiny{n/N}{h}  \combitiny{-Nh+(i+1)n/N}{n/N}}{t+Nh}  \mbox{ for } N\leq i \leq (r+1)N-1. 
\end{eqnarray*}
Then the partial fraction expansion of $F= F_0 ^{a-2rN} G_1\ldots G_{rN} H_N \ldots H_{(r+1)N-1}$ can be obtained my multiplying those of $F_0$, $G_i$ and $H_i$ using repeatedly the formula
\begin{equation} \label{eqprodelsples}
\frac1{(t+Nh)(t+Nh')^\ell } = \frac1{N^\ell (h'-h)^\ell (t+Nh)} - \sum_{i=1}^\ell \frac1{N^{\ell+1-i} (h'-h)^{\ell+1-i} (t+Nh')^i} 
\eneq
with $h\neq h'$. The denominator of $p_{j,h} $ comes both from this formula (and this contribution divides $(Nd_{n/N})^{a-j}$) and from the denominators of the coefficients in the partial fraction expansions of $F_0$, $G_i$, $H_i$ (which belong to $N^{-n/N}\Z$, so that $N^{an/N}$ accounts for this contribution). This concludes the proof of \eqref{eqdenompjh}.

\bigskip

On the other hand, bounding from above the coefficients of the 
partial fraction expansions of $F_0$, $G_i$, $H_i$  yields
$$|p_{j,h}| \leq n^{O(1)} N^{-an/N} 2^{an/N}\prod_{i=1}^{rN } \frac{(n+in/N)!}{(n/N)! (n+(i-1)n/N)!} \prod_{i=N}^{(r+1)N-1 } \frac{((i+1)n/N)!}{(n/N)! ( i n/N)!} $$
where $O(1)$ is a constant depending only on $a$, $r$, $N$ which can be made explicit (see \cite{Colmez} for details). Simplifying the products and using the bound $\frac{m!}{m_1!\ldots m_c!} \leq c^m$ valid when $m_1+\ldots+m_c = m$, one obtains
$$|p_{j,h}| \leq n^{O(1)} (2/N)^{an/N} \Big(\frac{((r+1)n)!}{n! (n/N)!^{rN}}\Big)^2 \leq  n^{O(1)} (2/N)^{an/N} (rN+1) ^{(2(r+1)n}.$$
This concludes the proof of Lemma \ref{lempjh}.

\bigskip

\Dem of  Proposition \ref{propformelin}:
Let $H(P)$ denote the exponential height of a polynomial $P\in\C[X]$, that is the maximum modulus of a coefficient of $P$. Recall that $P_j(z) = \sum_{h=0}^{n/N} p_{j,h} z^{Nh}$, $U(z)  = -\sum_{t=1}^n z^t \sum_{j=1}^a \sum_{h=0}^{\lfloor (t-1)/N\rfloor} \frac{p_{j,h}}{(t-Nh)^j}$
and $V(z)  = -\sum_{t=0}^{n-1} z^t \sum_{j=1}^a \sum_{h= \lceil (t+1)/N\rceil} ^{n/N}  \frac{p_{j,h}}{(Nh-t)^j}$.
Using Lemma \ref{lempjh} we see that these polynomials have coefficients in $\delta_n^{-1}\Z$ and height less than $H_n$ for some $H_n \leq 
\Big( (2/N)^{a/N} (rN+1)^{2r+2}\Big)^{n(1+o(1))}$. Now let $\Pti_{k,j} = z^{k-1} P_{k,j}$ for any $k$, $j$. Then the recurrence relation \eqref{eqdefpkj} yields
$$\Pti _{k,j} = z \Pti_{k-1,j}' - (k-2) \Pti_{k-1,j} - \Pti_{k-1,j+1}$$
where $\Pti_{k-1,j+1}=0$ if $j=a$, so that $\Pti_{k,j}$ is a polynomial of degree at most $n$, with coefficients in $\delta_n^{-1}\Z$ and height $H( \Pti_{k,j}) \leq (n+1)_{k-1} H_n$, by induction on $k$.

In the same way, letting $\Uti_k = z^{k-1} U_k$, Eq. \eqref{eqdefpz} yields
$$\Uti_k = z \Uti_{k-1}' - (k-2)\Uti_{k-1} -zQ_{k-1}$$
where $Q_{k-1} =\frac1{1-z} \Pti_{k-1,1}$. Provided $k\leq \dz-1$, Eq. \eqref{eqokenun} asserts that $P_{k-1,1}(1)=0$ so that $Q_{k-1}$ is a polynomial and $H(Q_{k-1})\leq n H(P_{k-1,1})\leq (n)_{k-1} H_n$. By induction on $k\leq \dz-1$, we deduce that $\Uti_k$ is   a polynomial of degree at most $n$, with coefficients in $\delta_n^{-1}\Z$ and height $H( \Uti_{k }) \leq k (n)_{k-1} H_n$. Now Eq. \eqref{eqpkzinfom}
reads 
$\Uti_k(z) = \sum_{\lam = 1}^N z^{\lam-1} U_{k,\lam}(z) $ with $U_{k,\lam} \in \Q[z^N, z^{-N}]$. If $k\leq \dz-1$  then $ U_{k,\lam}$ belongs to  $\Q[z^N]$, has degree at most $n$ (as a polynomial in $z$), coefficients in $\delta_n^{-1}\Z$ and height $H(U_{k,\lam}) \leq k (n)_{k-1} H_n$.

Proceeding in the same way, it is not difficult to prove that the same properties hold for $V_{k,\lam}$. Assertion $(i)$ follows at once, since $d_n = e^{n(1+o(1))}$.

\bigskip

\newcommand{\BAA}{{\mathfrak A}}

To prove $(ii)$, we recall that $\dz = -\deg F$ and write, as $|t| \to \infty$:
$$F(t) = \sum_{d = \dz} \frac{\BAA_d}{t^d} \mbox{ where  }
\BAA_d  = \sum_{j=1}^{ a }\sum_{h=0}^{n/N} (-Nh)^{d-j} \combi{d-1}{d-j} p_{j,h}$$
since $(t+Nh)^{-j} = \sum_{\ell = 0}^\infty \combitiny{\ell+j-1}{\ell} (-Nh)^\ell t^{-j-\ell}$ (see \cite[p. 1378]{FR}). Lemma \ref{lempjh} provides a positive real number $A_n \leq \Big( (2/N)^{a/N} (rN+1)^{2r+2}\Big)^{n(1+o(1))}$ such that $| \BAA_d | \leq (2n)^d A_n $ for any $d\geq \dz$. Then we have for any $t\in\Z$ such that $|t| \geq 2n+1$:
\begin{equation} \label{eqmajoF}
|F(t)| \leq A_n \sum_{d=\dz}^\infty (2n/t)^d \leq (2n+1)A_n (2n/t)^{\dz}.
\eneq
For any $z\in\C$ such that $|z| \leq 1$, and any $k\leq \dz-1$, we obtain
\begin{eqnarray*}
|S_0^{(k-1)}(z) | &=& | \sum_{t=(r+1)n+1}^\infty F(-t) (t-k+2)_{k-1} z^{t-k+1}| \leq (2n+1)A_n (2n)^{\dz} \sum_{t=(r+1)n+1}^\infty  t^{k-1-\dz}\\
&\leq&  (2n+1)A_n (2n)^{\dz}  \int_{(r+1)n}^{\infty}   t^{k-1-\dz} \dd t \leq  (2n+1)A_n  2 ^{\dz} n^k r^{k-\dz}.
\end{eqnarray*}
Moreover the same upper bound holds for $S_\infty(z) =  \sum_{t= r n+1}^\infty F( t) z^{-t}$ provided $|z| \geq 1$. Since 
$$S(z) = \sum_{\ell=1}^N \om^\ell \mu_\ell S_0(\om^\ell z) +  \om^\ell \nu_{N-\ell} S_\infty(\om^\ell z).$$
and $\dz = a(n/N+1)-2rn$, we obtain $| \delta_n   S^{(k-1)}(1)| \leq  \alpha ^{n(1+o(1)}$ for any $z\in \C$ such that $|z| = 1$, and any $k\leq \dz-1$; here the constant implied in $o(1)$ may depend on $k$ (but not on $n$). This concludes the proof of Proposition \ref{propformelin}.

\subsection{Siegel's linear independence criterion}\label{subsec36}

The proofs of all linear independence results in this paper rely on the following criterion, which is based on  Siegel's ideas (see for instance \cite[p. 81--82 and 215--216]{EMS}, \cite[\S 3]{Matala-Aho}  or  \cite[Proposition 4.1]{Marcovecchio}).

\begin{Prop} \label{propsiegel}
Let $\theta_1,\ldots,\theta_p$ be  real numbers, not all zero. Let $\tau>0$, and $(Q_n)$ be a sequence of real numbers with limit $+\infty$. Let $\calN$ be an infinite subset of $\N$, and for any $n\in\calN$ let $L^{(n)} = [\ell_{i,j}^{(n)}]_{1\leq i,j\leq p}$ be a matrix with integer coefficients and non-zero determinant, such that as $n\to\infty$ with $n\in\calN$:
$$\max_{1\leq i,j\leq p } | \ell_{i,j}^{(n)}| \leq Q_n^{1+o(1)}$$
$$
\mbox{ and } 
\max_{1\leq  j\leq p } | \ell_{1,j}^{(n)} \theta_1 + \ldots  + \ell_{p,j}^{(n)} \theta_p | \leq Q_n^{-\tau + o(1)}.$$
Then we have
$$\dim_\Q\Span_\Q(\theta_1,\ldots,\theta_p)\geq \tau+1.$$
\end{Prop}

In the proof of Theorem \ref{thprinc} we apply this proposition with $Q_n = \beta^n$ and $\tau = -\frac{\log \alpha}{\log\beta}$ (so that $Q_n^{-\tau} = \alpha^n$), where $\alpha $ and $\beta$ are defined in \S \ref{subsec32}; $\calN$ is the set of integer multiples of $N$.

\bigskip

Eventhough it is a classical result, let us recall the proof of Proposition \ref{propsiegel}. Let $d = \dim_\Q\Span_\Q(\theta_1,\ldots,\theta_p)$, and 
$F$ be a subspace of $\R^p$  defined over $\Q$, of dimension $d$,  which contains the point $(\theta_1,\ldots,\theta_p)$. 
  Let $n\in \calN$ be sufficiently large, and denote by $L_j^{(n)}$ the linear form $\ell_{1,j}^{(n)} X_1 + \ldots  + \ell_{p,j}^{(n)} X_p$ on $\R^p$. 
  Up to reordering 
$L_1^{(n)}$, \ldots, $L_p^{(n)}$, we may assume the restrictions of $L_1^{(n)}$, \ldots, $L_d^{(n)}$ to $F$ to be   linearly independent  linear forms on $F$. Denoting by $(u_1,\ldots,u_d)$ an $\R$-basis of $F$ consisting in vectors of $\Z^p$, the matrix $[L_j^{(n)} (u_t)]_{1\leq  j, t\leq d}$ has a non-zero integer determinant. Now $(\theta_1,\ldots,\theta_p)$ is a linear combination of $u_1$, \ldots, $u_d$; the same linear combination of the columns has coefficients less than $Q_n^{-\tau + o(1)}$ in absolute value. Therefore $ Q_n^{d-1-\tau + o(1)}$ is an upper bound on this non-zero integer determinant: this concludes the proof of Proposition \ref{propsiegel}.

\newcommand{\url}{\texttt}
\providecommand{\bysame}{\leavevmode ---\ }
\providecommand{\og}{``}
\providecommand{\fg}{''}
\providecommand{\smfandname}{\&}
\providecommand{\smfedsname}{\'eds.}
\providecommand{\smfedname}{\'ed.}
\providecommand{\smfmastersthesisname}{M\'emoire}
\providecommand{\smfphdthesisname}{Th\`ese}

\bigskip

S. Fischler, 
Laboratoire de Math\'ematiques d'Orsay, Univ. Paris-Sud, CNRS, Universit\'e Paris-Saclay, 91405 Orsay, France.

\end{document}